\documentclass[12pt]{article}
\usepackage{amsmath,amssymb,amsthm,eucal}
\usepackage{hyperref}
\usepackage{color}
\usepackage[normalem]{ulem}
\usepackage{enumerate}
\usepackage{microtype}
\usepackage{todonotes}

\font\tenrm=cmr10

\font\cmssl=cmss10 at 12 pt

\font\bigss=cmssdc10 scaled 2300

\font\cmsslll=cmss10 at 14 pt

\newcommand{\f}{\varphi}

\renewcommand{\o}{\omega}

\newcommand{\z}{\zeta}

\newcommand{\bC}{\mathbb{C}}
\newcommand{\bR}{\mathbb{R}}

\newcommand{\bN}{\mathbb{N}}

\newcommand\SO{\mathrm{SO}}

\newcommand{\p}{\partial}

\renewcommand{\square}{\kern1pt\vbox
 {\hrule height 0.6pt\hbox{\vrule width 0.6pt\hskip 3pt
 \vbox{\vskip 6pt}\hskip 3pt\vrule width 0.6pt}\hrule height0.6pt}
 \kern1pt}

\newcommand{\ra}{\rightarrow}

\newtheorem{Pb}{Problem}

\newtheorem{Th}{Theorem}
\newtheorem{Prop}[Th]{Proposition}
\newtheorem{Ex}[Th]{Example}

\newtheorem{Cor}[Th]{Corollary}
\newtheorem{Lem}[Th]{Lemma}
\newtheorem{Def}[Th]{Definition}

\theoremstyle{definition}

\newcommand{\bP}{\begin{Pb}\ \ }
\newcommand{\eP}{\end{Pb}}
\newcommand{\bt}{\begin{Th}\ \ }
\newcommand{\et}{\end{Th}}
\newcommand{\bp}{\begin{Prop}\ \ }
\newcommand{\ep}{\end{Prop}}
\newcommand{\bc}{\begin{Cor}\ \ }
\newcommand{\ec}{\end{Cor}}
\newcommand{\bl}{\begin{Lem}\ \ }
\newcommand{\el}{\end{Lem}}
\newcommand{\bd}{\begin{Def}\ \ }
\newcommand{\ed}{\end{Def}}
\newcommand{\pf}{\begin{proof}[{\it Proof:\ \ }]}
\newcommand{\epf}{\end{proof}}

\newcommand{\be}{\begin{equation}}
\newcommand{\ee}{\end{equation}}

\newcommand\re[1]{(\ref{#1})}
\newcommand{\arr}{\begin{array}{rlll}}
\newcommand{\ea}{\end{array}}
\newcommand{\bea}{\begin{eqnarray}}
\newcommand{\eea}{\end{eqnarray}}
\newcommand{\bean}{\begin{eqnarray*}}
\newcommand{\eean}{\end{eqnarray*}}

\catcode`@=11
\@addtoreset{equation}{section}
\catcode`@=12

\theoremstyle{remark}

 \newcommand{\zb}{\bar{z}}

\renewcommand{\z}{\zeta}
\newcommand{\zt}{\tilde{\zeta}}

\newcommand{\opn}{\operatorname} 
\newtheorem{Rem}[Th]{Remark} 
\newcommand{\br}{\begin{Rem}\ \ } 
\newcommand{\er}{\end{Rem}} 

\begin{document}
\rightline{}
\vskip 1.5 true cm
\begin{center}
{\bigss 
A class of cubic hypersurfaces and quaternionic\\[.5em] K\"ahler manifolds of co-homogeneity one}
\vskip 1.0 true cm
{\cmsslll V.\ Cort\'es, M.\ Dyckmanns, M.\ J\"ungling and D.\ Lindemann
} \\[3pt]
{\tenrm Department of Mathematics\\
and Center for Mathematical Physics\\
University of Hamburg\\
Bundesstra{\ss}e 55,
D-20146 Hamburg, Germany\\
vicente.cortes@uni-hamburg.de, malte.dyckmanns@uni-hamburg.de,\\ michel.juengling@googlemail.com, david.lindemann@uni-hamburg.de}\\[1em]
\vspace{2ex}
December 27, 2019\\
\end{center}
\vskip 1.0 true cm
\baselineskip=18pt
\begin{abstract}
\noindent
We classify all complete projective special real manifolds with reducible cubic potential, obtaining four series. 
For two of the series the manifolds are homogeneous, for the two others the respective automorphism group acts with
co-homogeneity one.
Complete projective special real manifolds give rise to complete quaternionic K\"ahler manifolds via the supergravity q-map, which is the composition of the supergravity c-map and r-map. We develop curvature formulas for manifolds in the image of the q-map. Applying the q-map to one of the above series of projective special real manifolds,
we obtain a series of complete quaternionic K\"ahler manifolds, which are shown to be inhomogeneous (of co-homogeneity one) 
based on our curvature formulas.\\[.5em] 
\end{abstract}
\tableofcontents

\section*{Introduction}
In this paper we are concerned with hypersurfaces $\mathcal{H} \subset \bR^{n+1}$ contained in the level set
$\{ h=1\}$ of a homogeneous cubic polynomial $h$. The hypersurface is equipped with the symmetric tensor field 
$g_{\mathcal{H}}$ on $\mathcal{H}$ induced by $-\frac{1}{3}\p^2h$. We require that $g_\mathcal{H}$ is a Riemannian metric. 
Then $(\mathcal{H},g_{\mathcal{H}})$ is called a \emph{projective special real manifold}, see Definition \ref{PSRDef}, $h$ is called its \emph{cubic potential} 
and $g_{\mathcal{H}}$ is called the \emph{projective special real metric}. 
The polynomials $h$ which admit such a hypersurface are called \emph{hyperbolic}, cf.\ Definition \ref{hyp:def}. 
Projective special real manifolds occur in the physics literature as the scalar manifolds of $5$-dimensional 
supergravity coupled to vector multiplets, see \cite{GST}. These manifolds are related to \emph{projective special 
K\"ahler manifolds} \cite{F,ACD}, reviewed in Section \ref{ConicalSect}, by a construction known as the \emph{r-map} \cite{DV}, which is induced by the dimensional reduction of the supergravity theory from $5$ to $4$ space-time dimensions. 

Similarly, projective special K\"ahler manifolds
are related to \emph{quaternionic K\"ahler manifolds} of negative scalar curvature 
by the \emph{$c$-map}, which is induced by dimensional reduction to $3$ dimensions \cite{FS}. 
As shown in \cite{ACM,ACDM}, the quaternionic K\"ahler property of the c-map metric can be proven by 
showing that it is part of a one-parameter family of metrics arising from an 
indefinite version of Haydys' HK/QK-correspondence \cite{Ha,Hi}. The key fact is that the 
cotangent bundle of any conical affine special K\"ahler manifold admits an indefinite hyper-K\"ahler metric 
(known as the rigid c-map metric) and 
a circle action of the type required for the correspondence to yield a family of positive definite quaternionic K\"ahler metrics. 
These constructions are recovered in \cite{MS} based on Swann's twist construction and elementary deformation. See also 
\cite{APP}, where it was first shown that the c-map metric and the rigid c-map metric are related by the 
QK/HK-correspondence. 

It is known \cite{CHM} that the r- and c-map preserve the completeness of the underlying Riemannian metrics.\footnote{The same is true for the 
\emph{generalized} r-map \cite{CHM}, where the cubic polynomial is replaced by a more general homogeneous function $h$. However, the resulting K\"ahler manifolds 
are in general no longer projective special K\"ahler and therefore not interesting for our present purposes.} It follows that the same is true
for their composition, the \emph{q-map}. In this way the study of the completeness of quaternionic K\"ahler manifolds obtained 
by the q-map is reduced to the study of the completeness of the initial projective special real manifold. Complete
projective special real manifolds are characterized by the following theorem, to be used later on. 
\begin{Th}[{\cite[Thm.\ 2.5]{CNS}}]\label{CNS:Thm}
A projective special real manifold $\mathcal{H}\subset \bR^{n+1}$ is complete with respect to the metric $g_{\mathcal{H}}$ if and only if $\mathcal{H}$ is closed as a subset of 
$\mathbb{R}^{n+1}$.
\end{Th}

It follows from Theorem \ref{CNS:Thm} that the classification of complete projective special real manifolds is equivalent to the solution
of the following two problems:
\begin{enumerate}[(i)]
\item Classification of all hyperbolic homogeneous cubic polynomials $h$, up to linear transformations. 
\item For each such polynomial determine all locally strictly convex components of the level set $\{ h=1 \}$, up to 
linear transformations. 
\end{enumerate}
While it is certainly possible to solve these problems in low dimensions, see \cite{CDL} for the solution up to polynomials in $3$ variables, we do not expect a simple solution 
valid in all dimensions. A very rough idea about problem (i) is obtained by observing that the dimension of the space 
of homogeneous cubic polynomials grows cubically whereas the dimension of the general linear group grows only 
quadratically with the number of variables. Notice that the hyperbolic polynomials form an open subset in the space of 
homogeneous cubic polynomials in a given number of variables. An interesting class of projective special real manifolds
is provided by considering those with reducible cubic potentials $h$, that is $h$ is a product of polynomials of lower degree. 
The motivation to consider this class is that the polynomial $h$ is preserved by a large group of linear transformations. 
In virtue of the general results about symmetries in Appendix \ref{appendix_A}, this will eventually give rise to quaternionic 
K\"ahler manifolds with a large but not always transitive group of isometries.  
Applying the q-map to the complete manifolds in this class we obtain a class of complete quaternionic K\"ahler manifolds, as 
follows from the general result \cite[Thm.\ 6]{CHM}. In this way one obtains, in particular, the series of symmetric spaces 
\begin{equation}
 \label{so:Eq} \frac{\SO_0(4,m)}{\SO(4)\times \SO(m)},\quad m\ge 3,
\end{equation}
as well as the series of homogeneous non-symmetric spaces $\mathcal{T}(p)$, $p\ge 1$, of rank $3$, 
see \cite{DV,C}. One of the results of this paper is that one also obtains a series of complete quaternionic K\"ahler manifolds that are not
locally homogeneous, see Theorem \ref{nonlochomQKseries}. In fact, we show that there are precisely four series of complete projective special real manifolds 
with reducible cubic potential.
More precisely, by solving the above problems (i) and (ii) under the assumption that $h$ is reducible we will 
obtain the following result.
\bt
\label{Thm2}
Every complete projective special real manifold $\mathcal{H} \subset \{ h=1\}\subset \bR^{n+1}$ of dimension 
$n\ge 2$ for which $h$ is reducible is linearly equivalent to exactly one of the following complete 
projective special real manifolds:
\begin{itemize}
\item[a)] $\{ x_{n+1}(\sum_{i=1}^{n-1} x_i^2-x_n^2)=1 ,\quad x_{n+1}<0, x_n>0\}$, 
\item[b)] $\{ (x_1+ x_{n+1})(\sum_{i=1}^{n} x_i^2-x_{n+1}^2)=1 ,\quad x_1+x_{n+1}<0\}$,
\item[c)] $\{ x_1(\sum_{i=1}^{n} x_i^2-x_{n+1}^2)=1 ,\quad x_1<0, x_{n+1}>0\}$,
\item[d)] $\{ x_1(x_1^2-\sum_{i=2}^{n+1} x_i^2)=1 ,\quad x_1>0\}$.
\end{itemize}
\et 

Notice that in the case $n=2$ the result follows from \cite[Thm.\ 1]{CDL} and that the above list is also valid in the case 
$n=1$ but then the curves a) and b) are linearly equivalent, as well
as c) and d), see \cite[Cor.\ 4]{CHM}. 

Under the q-map the series a) with $n\ge 1$ corresponds to the series \re{so:Eq} of symmetric quaternionic K\"ahler manifolds with $m=n+2$. 
Similarly, b) corresponds to the series $\mathcal{T}(p)$ of homogeneous quaternionic K\"ahler manifolds with $p=n-1\ge 0$, where only
the first member $\mathcal{T}(0)= \frac{\SO_0(4,3)}{\SO(4)\times \SO(3)}$ of the series is symmetric. The quaternionic K\"ahler manifolds 
obtained from the series c) and d) admit a Lie group acting isometrically with co-homogeneity one. For d) we will prove the following 
stronger result.

\bt \label{inh:Thm} The quaternionic K\"ahler manifolds associated with the projective special
real manifolds $\{ (x_1,\ldots,x_{n+1})\in \bR^{n+1}\mid x_1(x_1^2-\sum_{i=2}^{n+1} x_i^2)=1 ,\quad x_1>0\}$, $n\ge1$, are complete of 
negative scalar curvature and the isometry group acts with co-homogeneity one. 
\et 

Notice that the theorem provides examples of complete quaternionic K\"ahler manifolds in all dimensions $\ge 12$ 
for which the isometry group acts with co-homogeneity one, see \cite{DS,PV} for results excluding the existence of such manifolds in the case of positive scalar curvature, and 
\cite{P} for some examples of negative scalar curvature in dimension $4$. Recall also that the quaternionic hyperbolic space admits deformations by complete quaternionic K\"ahler manifolds \cite{L}, but the isometry
groups of these are not known.

The claim that the quaternionic K\"ahler manifolds in Theorem \ref{inh:Thm}, and similarly the ones obtained from the series c) in Theorem \ref{Thm2}, 
admit a \emph{subgroup} of the 
isometry group acting with an orbit of codimension one follows from the fact that the automorphism group of the initial projective special
real manifolds acts with an orbit of codimension one.
In fact, the orthogonal group $\mathrm{O}(n)$ in the variables $x_2,\ldots,x_{n+1}$ acts by automorphisms of the projective special real manifold. Moreover, every automorphism of a projective special
real manifold extends to an isometry of the corresponding quaternionic K\"ahler manifold under the q-map. In addition, 
the r-map as well as the c-map each produce a freely acting additional solvable Lie group of automorphisms, see \cite{DV,DVV,CHM}. 
The dimensions of the latter solvable groups coincide with the number of extra dimensions created by the 
r- and c-map, respectively. Therefore the co-homogeneity does not increase under these constructions, see Appendix \ref{appendix_A} for details.

The main difficulty is to prove that the quaternionic K\"ahler manifolds of Theorem \ref{inh:Thm} are not of co-homogeneity zero, this is the content of Theorem \ref{nonlochomQKseries}. 
The proof proceeds by computing the point-wise norm of the curvature tensor and showing that for each of these manifolds it is a non-constant rational function 
depending only on one coordinate $x$ out of a system of $4n+8$ global coordinates. It relies on general
curvature formulas for quaternionic K\"ahler manifolds obtained by the q-map, which constitute another important result of this paper, see Theorem \ref{thQuarticTensorQMap} and Corollary \ref{Rsquarenorm}. Incidentally, we expect that the isometry groups of the quaternionic 
K\"ahler manifolds corresponding to the remaining series c) in Theorem \ref{Thm2} do likewise have co-homogeneity precisely one. 
The corresponding curvature calculations are more involved in that case. 


\subsubsection*{Acknowledgements}
This work was
partly supported by the German Science Foundation (DFG) under the Research Training Group 1670 ``Mathematics inspired by String Theory". 
V.C. would like to thank Thomas Mohaupt for many discussions about isometries of r- and c-map spaces
in the framework of our joint research projects.

\section{Classification of complete projective special real manifolds with redu\-cible cubic potential}
In this section we will classify all complete projective special real manifolds with reducible cubic potential up to linear transformations. 
After giving some basic definitions we will first classify up to equivalence all non-degenerate reducible homogeneous cubic polynomials in Section
\ref{redSec} and among these all hyperbolic ones in Section \ref{hypSec}. In the same section we determine, for each of the resulting hyperbolic polynomials $h$, those connected components (up to
linear transformations) of the level sets $\{ h=1\}$ which contain a hyperbolic point, see Definition \ref{hyp:def}. In particular we 
determine all such components which are locally strictly convex or, equivalently, consist solely of hyperbolic points. 
As a consequence of Theorem \ref{CNS:Thm} these components give precisely all complete projective special real manifolds 
with reducible cubic potential (up to linear transformations). 

\bd \label{hyp:def}
Let $h:\mathbb{R}^{n+1}\to\mathbb{R}$ be a homogeneous cubic polynomial. 
\begin{enumerate}
\item The polynomial 
$h$ is called {\cmssl non-degenerate} if there exists $p\in\mathbb{R}^{n+1}$, such that $\det \partial^2 h_p\ne 0$. 
\item The polynomial $h$ is called {\cmssl hyperbolic} if there exists a {\cmssl hyperbolic point} $p\in \bR^{n+1}$, that is a point such that $h(p)>0$ and $\partial^2 h_p$ is of signature $(1,n)$. 
\end{enumerate}
Two homogeneous cubic polynomials are called {\cmssl equivalent} if they are related by a linear transformation. 
\ed
Notice that the notions of non-degeneracy and hyperbolicity are invariant under linear transformations and that $\det \partial^2 h_p\ne 0$ implies $h(p)\neq 0$ if the 
tensor $\partial^2 h_p$ is non-degenerate on the hyperplane $\ker dh_p$.

\bd \label{PSRDef}
A hypersurface $\mathcal{H}\subset\mathbb{R}^{n+1}$ is called a {\cmssl projective special real manifold} if there exists a homogeneous cubic polynomial $h:\mathbb{R}^{n+1}\to\mathbb{R}$, such that
\begin{enumerate}[(i)]
\item $\mathcal{H}\subset\{x\in\mathbb{R}^{n+1}\ |\ h(x)=1\}$ and
\item $g_{\mathcal{H}}:=-\frac13 \partial^2h|_{T\mathcal{H}\times T\mathcal{H}}>0$.
\end{enumerate}
The hypersurface $\mathcal{H}\subset\mathbb{R}^{n+1}$ is endowed with the Riemannian metric $g_{\mathcal{H}}$ which is called
the {\cmssl projective special real metric}\footnote{For practical reasons, we prefer to compute $-\frac12\partial^2h$ instead of $-\frac{1}{3}\partial^2h$ below.} or {\cmssl centroaffine metric}, see \cite{CNS} for an explanation of this terminology. 
Two projective special real manifolds are called {\cmssl isomorphic} if there is a linear transformation inducing a bijection between them. 
\ed
\br Using the homogeneity of $h$ one sees that for every projective special real manifold $\mathcal{H}$ the symmetric tensor 
$\partial^2 h_p$ is of signature $(1,n)$ for all $p\in \mathcal{H}$ and that 
$\mathcal{H}$ is perpendicular to the position vector $p$ with respect to $\partial^2 h_p$. 
In fact, $\partial^2h_p(p , p) = 6 h(p) = 6$ and 
$\partial^2h_p (p , v) = 3 dh_pv=0$ for all $p\in \mathcal{H}$, $v\in T_p\mathcal{H}$. In particular, $h$ is hyperbolic. 
Notice also
that a linear transformation mapping a projective special real manifold $\mathcal{H} \subset \bR^{n+1}$ to another 
projective special real manifold $\mathcal{H}' \subset \bR^{n+1}$ is automatically an isometry with respect to the centroaffine 
metrics. In particular, isomorphic projective special real manifolds are isometric.
\er

 In order to avoid special cases in low dimensions, and since the case $n\le 2$ has already been studied \cite{CDL}, we will always assume that $n\ge 3$ in the following classifications. 

\subsection{Classification of non-degenerate reducible po\-lynomials}
\label{redSec}

For $m\in \bN$ and $k\in\{0,\ldots,m\}$, we introduce the following quadratic polynomials on $\mathbb{R}^{m}$:
\begin{equation*}
Q^m_k:=\sum\limits_{i=1}^{k}x_i^2-\sum\limits_{i=k+1}^{m}x_i^2. 
\end{equation*}

\bp \label{nondegProp} 
Any non-degenerate reducible homogeneous cubic polynomial $h$ on $\mathbb{R}^{n+1}$, $n \geq 3$, is equivalent to 
precisely one of the following:
\begin{enumerate}[I)]
\item $x_{n+1} Q^n_k$, $\frac{n}{2}\leq k\leq n$,
\item $x_1 Q^{n+1}_k$, $1\leq k\leq n+1$,
\item $(x_1+x_{n+1}) Q^{n+1}_k$, $\frac{n+1}{2}\leq k\leq n$.
\end{enumerate}
\ep

\begin{proof}
Let $h=LQ$ be a non-zero reducible cubic polynomial on $\mathbb{R}^{n+1}$, where $L$ is a linear and $Q$ a quadratic factor. 
Up to a linear transformation, we can assume that $Q=Q^m_k$, $1\leq m\leq n+1$, $\frac{m}{2}\leq k \leq m$. 
In the following, let
\begin{equation*}
L :=\sum_{j=1}^{n+1} a_jx_j.
\end{equation*}
Next we examine for which choices of $Q^m_k$ and $L$ the polynomial $h=LQ^m_k$ is non-degenerate.
Notice that $m=n$ or $m=n+1$, since otherwise $0\neq \ker dL \cap \ker \p^2 Q \subset \ker \p^2 h_p$ for all $p\in \bR^{n+1}$.

In the case $m=n$ the non-degeneracy of $h$ clearly implies that $a_{n+1}\neq 0$ and without loss of generality we can
assume that $L=x_{n+1}$. We compute 
\begin{equation*}
\partial^2h=2\left(\begin{matrix}
x_{n+1} & & & & & & x_1 \\
 & \ddots & & & & & \vdots \\
 & & x_{n+1} & & & & x_k \\
 & & & -x_{n+1} & & & -x_{k+1} \\
 & & & & \ddots & & \vdots \\
 & & & & & -x_{n+1} & -x_n \\
x_1 & \hdots & x_k & -x_{k+1} & \hdots & -x_n & 0
\end{matrix}\right), 
\end{equation*}
where the remaining entries are zero. The determinant is given by
\[ \det \partial^2h = 2^{n+1}(-1)^{n-k+1}x_{n+1}^{n-2}h,\]
which shows that $h=x_{n+1}Q^n_k$ is non-degenerate for all $\frac{n}{2}\leq k\leq n$. These are precisely the polynomials listed in I).

It remains to check the case $m=n+1$, that is, $h=LQ^{n+1}_k$, $\frac{n+1}{2}\leq k \leq n+1$. Using the transitive action of the pseudo-orthogonal group
of the quadratic form $Q^{n+1}_k$ on each pseudo-sphere and on the cone of non-zero light-like vectors we can assume up to a positive rescaling that 
$L=x_1$ ($L$ space-like), $L=x_{n+1}$ ($L$ time-like), or $L=x_1+x_{n+1}$ ($L$ light-like), where the latter two cases need only to be considered for $k\le n$.
Since $x_{n+1}$ is space-like with respect to $-Q^{n+1}_k$ for $\frac{n+1}{2} \le k\le n$ and $-Q^{n+1}_k$ is equivalent to $Q_{n+1-k}^{n+1}$, $1\le n+1-k \le \frac{n+1}{2}$, 
we are left with the two cases II) and III).

In case II), $h=x_1Q^{n+1}_k$ with $1\leq k\leq n+1$ and
\begin{equation*}
\partial^2h=2\left(\begin{matrix}
3x_1& x_2 & \hdots & x_k & -x_{k+1} & \hdots & -x_{n+1} \\
x_2 & x_1 & & & & & \\
\vdots & & \ddots & & & & \\
x_k & & & x_1 & & & \\
 -x_{k+1} & & & & -x_1 & & \\
\vdots & & & & & \ddots & \\
-x_{n+1} & & & & & & -x_1 \\
\end{matrix}\right).
\end{equation*}
We obtain
\begin{equation*}
\det\partial^2h=(-1)^{n+1-k}2^{n+1}x_1^{n-2}(4x_1^3-h), 
\end{equation*}
which, for all $1\leq k\leq n+1$, is not the zero polynomial. Hence, all polynomials listed in II) are non-degenerate.

In case III), that is $h=(x_1+x_{n+1})Q^{n+1}_k$, $\frac{n+1}{2}\leq k\leq n$, it is convenient to change the coordinates the following way:
\begin{align*}
x_1+x_{n+1}&=\xi,\\
x_1-x_{n+1}&=\eta.
\end{align*}
$h$ is now of the form
\begin{equation*}
h=\xi\left(\xi\eta+\sum_{i=2}^{k}x_i^2-\sum_{i=k+1}^n x_i^2\right) .
\end{equation*}
In the coordinates $(\xi , \eta , x_2,\ldots , x_n)$ we have 
\begin{equation*}
\partial^2h=
2\left(
\begin{array}{cccccccc}
\eta & \xi & x_2 & \hdots & x_k & -x_{k+1}&\hdots & -x_n\\ 
\xi & 	0		& 	 		&				&			&					&				&			\\ 
x_2 &				& \xi&				&			&					&				&			\\
\vdots&			&			& \ddots & 			&					&				&			\\
x_k	 &			&			&				& \xi&					&				&			\\
-x_{k+1}&		&			&				&			& -\xi		&				&			\\
\vdots &			&			&				&			&					& \ddots	&			\\
-x_n &			&			&				&			&					&				&-\xi\\
\end{array}
\right).
\end{equation*}
It is now easy to see that 
\begin{equation*}
\det \partial^2h = (-1)^{n+1-k}\xi^{n+1}.
\end{equation*}
We conclude that all polynomials considered in III) are non-degenerate.
\end{proof}

\subsection{Classification of hyperbolic reducible polynomials and 
\label{hypSec}
complete projective special real manifolds}
Let $h:\mathbb{R}^{n+1}\to \mathbb{R}$ be a hyperbolic homogeneous cubic polynomial. We consider the open subset $\mathcal{H}(h)$ of the hypersurface $\{h=1\}$ 
consisting of the hyperbolic points of $h$: 
\begin{equation*}
\mathcal{H}(h)=\{p\in\mathbb{R}^{n+1}\ |\ h(p)=1,-\partial^2h_p\text{ has Lorentzian signature}\; (n,1)\}.
\end{equation*}
\bp
Let $h:\mathbb{R}^{n+1}\to \mathbb{R}$, $n\geq 3$, be a reducible hyperbolic homogeneous cubic polynomial and let $(x_1,\hdots,x_{n+1})$ denote the standard coordinates of $\mathbb{R}^{n+1}$. Then $h$ is equivalent to one of the following polynomials and the corresponding hypersurface $\mathcal{H}(h)$ endowed with the Riemannian metric $-\frac{1}{2} \partial^2h|_{T\mathcal{H}(h)\times T\mathcal{H}(h)}$ has the following properties:
\begin{enumerate}[a)]
\item $h=x_{n+1}\left( \sum\limits_{i=1}^{n-1} x_i^2-x_n^2 \right)$, $\mathcal{H}(h)=\{h=1,\ x_{n+1}<0\}$ has two connected components, both closed and isomorphic.
\item $h=(x_1+x_{n+1})\left( \sum\limits_{i=1}^n x_i^2 -x_{n+1}^2 \right)$, $\mathcal{H}(h)=\{h=1,\ x_1+x_{n+1}<0\}$ has one connected component and it is closed.
\item $h=x_1\left( \sum\limits_{i=1}^n x_i^2-x_{n+1}^2\right)$, $\mathcal{H}(h)=\{h=1,\ x_1<0\}$ has two connected components, both closed and isomorphic.
\item $h=x_1\left( x_1^2 -\sum\limits_{i=2}^{n+1} x_i^2 \right)$, $\mathcal{H}(h)=\{h=1,\ x_1>0\}$ has one connected component and it is closed. 
\item $h=x_1\left( x_1^2+x_2^2-\sum\limits_{i=3}^{n+1} x_i^2 \right)$, $\mathcal{H}(h)=\{h=1\}\cap \{\frac{1}{\sqrt[3]{4}}>x_1>0\}$ has two connected components. 
They are isomorphic and not closed.
\end{enumerate}
In particular, the closed connected components of the respective $\mathcal{H}(h)$ are complete projective special real manifolds.
\begin{proof} In Proposition \ref{nondegProp} we have listed all non-degenerate cubic homogeneous polynomials up to equivalence. It remains to determine which ones are hyperbolic and to analyse the properties of the connected components of $\mathcal{H}(h)$. In the following we treat each of the cases I-III) of Proposition \ref{nondegProp}. 

{\bf I)} 
Recall that the family I) of Proposition \ref{nondegProp} contains the polynomials $h=x_{n+1}Q^n_k$, $\frac{n}{2}\leq k\leq n$, with
\begin{equation*}
-\frac{1}{2}\partial^2h=-x_{n+1}\left( \sum_{i=1}^k dx_i^2-\sum_{i=k+1}^n dx_i^2\right) 
-2\left( \sum_{i=1}^k x_idx_i-\sum_{i=k+1}^n x_idx_i\right) dx_{n+1}.
\end{equation*}
To check that a point $p\in\mathbb{R}^{n+1}$ is hyperbolic it suffices to construct an orthogonal basis of $T_p\mathbb{R}^{n+1}$ with respect to $-\frac{1}{2}\partial^2 h$ and to check that the Gram matrix has Lorentzian signature. Note that the vectors $\{\partial_{x_1},\ldots,\partial_{x_n}\}$ are orthogonal at each point:
\begin{equation*}
-\frac{1}{2}\partial^2h(\partial_{x_i},\partial_{x_j})=\left\{
\begin{array}{rc}
-\delta_i^jx_{n+1}, & 1\leq i,j\leq k,\\
\delta_i^jx_{n+1}, & k+1\leq i,j\leq n,\\
0, & \text{otherwise}.
\end{array} \right.
\end{equation*}
Now the restrictions $n\geq 3$, $k\geq \frac{n}{2}$, allow us to limit the possibility of hyperbolic points to the cases $k=n-1$ and $k=n$ and we obtain the requirement $x_{n+1}<0$. Other\-wise we would have at least two time-like vectors in an orthogonal basis of the form $(v,\partial_{x_1},\ldots,\partial_{x_n})$. For $v=\sum_{i=1}^{n+1}v_i\partial_{x_i}$ to be orthogonal to $\partial_{x_i}$ for all $1\leq i\leq n$ it has to fulfil
\begin{equation*}
x_{n+1}v_i+x_iv_{n+1}=0\ \ \forall 1\leq i\leq n.
\end{equation*}
Hence, $v_i=-\frac{x_iv_{n+1}}{x_{n+1}}$ for $1\leq i\leq n$ and $v=v_{n+1}\left( -\sum_{i=1}^n\frac{x_i}{x_{n+1}}\partial_{x_i} + \partial_{x_{n+1}} \right)$. Since $x_{n+1}<0$, we might choose $v=\sum_{i=1}^nx_i\partial_{x_i} - x_{n+1}\partial_{x_{n+1}}$ and obtain
\begin{equation*}
-\frac{1}{2}\partial^2h(v,v)=x_{n+1}\left( \sum_{i=1}^k x_i^2-\sum_{i=k+1}^n x_i^2 \right)=h.
\end{equation*}
Hyperbolic points need to fulfil $h(p)>0$ by definition, which implies $-\frac{1}{2}\partial^2h(v,v)>0$. Hence, $h=x_{n+1}Q^n_k$, $\frac{n}{2}\leq k\leq n$, is hyperbolic if and only if $k=n-1$, that is $h=x_{n+1}\left( \sum\limits_{i=1}^{n-1} x_i^2-x_n^2 \right)$ is the polynomial \textcolor{red}{$a)$} of this proposition. The hypersurface $\mathcal{H}(h)$ consists of the connected components
\begin{equation*}
\mathcal{H}_1:=\left\{ (x_1,\ldots,x_{n+1})\in\mathbb{R}^{n+1}\mid h(x_1,\ldots,x_{n+1})=1,\ x_n<0,x_{n+1}<0\right\}
\end{equation*}
and
\begin{equation*}
\mathcal{H}_2:=\left\{ (x_1,\ldots,x_{n+1})\in\mathbb{R}^{n+1}\mid h(x_1,\ldots,x_{n+1})=1,\ x_n>0,x_{n+1}<0\right\}.
\end{equation*}
One can easily verify that $\mathcal{H}_1$ and $\mathcal{H}_2$ are both closed in $\mathbb{R}^{n+1}$ and related by the involution $(x_1,\ldots,x_n,x_{n+1})\mapsto (x_1,\ldots,-x_n,x_{n+1})$. 

\paragraph*{II)}
The family II) of Proposition \ref{nondegProp} contains polynomials of the form $h=x_1Q^{n+1}_k$, $1\leq k \leq n+1$. We will construct an orthogonal basis for each $p\in\{h>0\}$, $p=(x_1,\ldots,x_{n+1})$, with respect to
\begin{equation*}
-\frac{1}{2}\partial^2h=x_1\left( -3dx_1^2-\sum_{i=2}^k dx_i^2 +\sum_{i=k+1}^{n+1} dx_i^2\right) -2dx_1\left( \sum_{i=2}^k x_idx_i-\sum_{i=k+1}^{n+1}x_i dx_i\right) .
\end{equation*}
We define 
\[ v = x_1\p_{x_1} -\sum_{i=2}^{n+1}x_i\p_{x_i}.\] 
Then one can check, for $x_1\neq 0$, that $(v,\partial_{x_2},\ldots,\partial_{x_{n+1}})$ is an orthogonal basis with respect to $-\frac{1}{2}\partial^2h$ and that 
\begin{equation*}
-\frac{1}{2}\partial^2h(v,v) = -4x_1^3+h.
\end{equation*} 
Thus, the possible values for $k$ that do not exclude the possibility for $h$ to be hyperbolic, the respective requirements for the possibly hyperbolic points, and the corresponding polynomials are (recall $n\geq 3$):
\begin{equation*}
\begin{array}{lllll}
A) &k=1, & x_1>0, & h<4x_1^3\ \, \left(-\frac{1}{2}\partial^2h(v,v)<0\right); & h=x_1\left(x_1^2-\sum_{i=2}^{n+1}x_i^2\right),\\ \\
B) &k=2, & x_1>0, & h>4x_1^3\ \, \left(-\frac{1}{2}\partial^2h(v,v)>0\right); & h=x_1\left( x_1^2+x_2^2-\sum_{i=3}^{n+1}x_i^2\right),\\ \\
C) &k=n, & x_1<0, & h>4x_1^3\ \, \left(-\frac{1}{2}\partial^2h(v,v)>0\right); & h=x_1\left(\sum_{i=1}^nx_i^2-x_{n+1}^2\right),\\ \\
D) &k=n+1, & x_1<0, & h<4x_1^3\ \, \left(-\frac{1}{2}\partial^2h(v,v)<0\right); & h=x_1\left(\sum_{i=1}^{n+1}x_i^2\right).
\end{array}
\end{equation*}
The polynomials in A), B), and C) are, in fact, hyperbolic, as seen by specifying a hyperbolic point: 
\begin{equation*}
\begin{array}{lll}
A) & p_A=(1,0,\ldots,0), & h(p_A)=1,\\ \\
B) & p_B=(1,2,0,\ldots,0), & h(p_B)=5,\\ \\
C) & p_C=(-1,0,\ldots,0,2), & h(p_C)=3.
\end{array}
\end{equation*}
These three series of polynomials are, corresponding to the above order $A)$, $B)$, and $C)$, the first three cases $d)$, $e)$, and $c)$ of this proposition.
The polynomials in $D)$ are not hyperbolic, since the specified conditions are not compatible with $h>0$. 
We will now describe the sets $\mathcal{H}(h)$.

In case $A)$, the set of hyperbolic points of $\bR^{n+1}$ with respect to $h$ was described by the inequalities
$x_1>0$ and $h<4x_1^3$. The second inequality follows from the first since $Q_1^{n+1}\le x_1^2$. This shows 
that $\mathcal{H}(h)=\{h=1,\ x_1>0\}$, which has one connected component. To see this consider for fixed $u=(x_2,\ldots,x_{n+1})\in\mathbb{R}^{n}$ the function 
\begin{equation*}
(\rho , \infty) \to\mathbb{R},\ x_1\mapsto h(x_1,u),
\end{equation*}
where $\rho= |u|$ 
and notice that it is a strictly monotonously increasing diffeomorphism onto $(0,\infty )$. In particular, for all $u\in \bR^n$ there is a unique $x_1(u)\in (\rho ,\infty )$ such that $h(x_1(u),u)=1$. We obtain a bijection 
\[ \bR^n \ra \mathcal{H}(h),\quad u\mapsto (x_1(u),u),\] 
which is a diffeomorphism by the implicit function theorem. In particular, $\mathcal{H}(h)$ is connected. This implies that it is a connected component of $\{h=1\}$ and, thus, closed in $\bR^{n+1}$.

In case $B)$, the requirement for hyperbolicity on $\{h=x_1\left( x_1^2+x_2^2-\sum_{i=3}^{n+1}x_i^2\right)=1\}$ is $\frac{1}{\sqrt[3]{4}}>x_1>0$, which implies 
 $x_2\ne 0$. Observe that
\begin{equation*}
h=1\ \Leftrightarrow\ x_2^2=\frac{1}{x_1}(1-x_1^3)+\sum_{i=3}^{n+1}x_i^2.
\end{equation*}
Hence, $\mathcal{H}(h)=\{h=1\}\cap \{\frac{1}{\sqrt[3]{4}}>x_1>0\}$ has two connected components, namely $\{h=1\}\cap \{\frac{1}{\sqrt[3]{4}}>x_1>0\}\cap\{x_2>0\}$ and $\{h=1\}\cap \{\frac{1}{\sqrt[3]{4}}>x_1>0\}\cap \{x_2<0\}$. They are related by the involution $x_2\mapsto -x_2$, which preserves the polynomial $h$. 
The two components of $\mathcal{H}(h)$ are not closed in $\mathbb{R}^{n+1}$, since its boundary is given by
\begin{equation*}\partial\mathcal{H}(h) =\left\{h=1,\, x_1=\frac{1}{\sqrt[3]{4} },\, \det \partial^2h =0\right\} 
=\left\{h=1,\, x_1=\frac{1}{\sqrt[3]{4} }\right\} =\left\{x_2^2-\sum_{i=3}^{n+1}x_i^2= \frac{3}{4^{\frac23}}
\right
\}\!\:\!.
\end{equation*}

In case $C)$, the requirement $x_1<0$ automatically implies the second requirement $h>4x_1^3$ on $\{h=x_1\left(\sum_{i=1}^nx_i^2-x_{n+1}^2\right)=1\}$ and, hence, $\mathcal{H}(h)=\{h=1,\ x_1<0\}$. Note that $\{h=1\}\cap \{x_1=0\}=\emptyset$ implies that the connected components of $\mathcal{H}(h)$ are also connected components
 of $\{h=1\}$, and thus are closed. $x_1<0$ and $h=x_1\left(\sum_{i=1}^nx_i^2-x_{n+1}^2\right)=1$ implies $\sum_{i=1}^nx_i^2-x_{n+1}^2<0$, which implies $x_{n+1}\ne 0$. Hence, the connected components of $\mathcal{H}(h)$ are given by the two graphs $\{h=1,\ x_1<0,\ x_{n+1}>0\}$ and $\{h=1,\ x_1<0,\ x_{n+1}<0\}$. They are related by the involution $x_{n+1}\mapsto-x_{n+1}$.

\paragraph*{III)}
Recall that each $h=(x_1+x_{n+1})Q^{n+1}_k$ contained in family III) of Proposition \ref{nondegProp} is equivalent to $h=\xi\left(\xi\eta+\sum_{i=2}^k x_i^2-\sum_{i=k+1}^n x_i^2\right)$. In these coordinates
\begin{equation*}
\begin{aligned}
-\frac{1}{2}\partial^2 h =\ &-\eta d\xi^2 -2\xi d\eta d\xi +\left( -2\sum_{i=2}^k x_i dx_i + 2 \sum_{i=k+1}^n x_i dx_i\right)d\xi\\
& +\xi \left( -\sum_{i=2}^k dx_i^2 +\sum_{i=k+1}^n dx_i^2 \right).
\end{aligned}
\end{equation*}
The set $\left\{h=\xi\left(\xi\eta+\sum_{i=2}^k x_i^2-\sum_{i=k+1}^n x_i^2\right)=1\right\}$ consists of exactly two connected components: \begin{equation*}
\mathcal{H}_1:= \left\{(\xi,\eta,x_2,\ldots,x_n)\in\mathbb{R}^{n+1}\left|\ \eta = \frac{1-\xi\left( \sum_{i=2}^{k}x_i^2-\sum_{i=k+1}^n x_i^2 \right)}{\xi^2},\ \xi>0\right\}\right.
\end{equation*}
and
\begin{equation*}
\mathcal{H}_2:= \left\{(\xi,\eta,x_2,\ldots,x_n)\in\mathbb{R}^{n+1}\left|\ \eta = \frac{1-\xi\left( \sum_{i=2}^{k}x_i^2-\sum_{i=k+1}^n x_i^2 \right)}{\xi^2},\ \xi<0\right\}\right. .
\end{equation*}
In order to determine which of the polynomials in this family are hyperbolic, we will pull back $-\frac{1}{2}\partial^2h$ to $\mathcal{H}_1$ and $\mathcal{H}_2$, respectively. We will use that $h$ is hyperbolic if and only if the pullback is Riemannian at least at one point contained in 
$\{h=1\}$. We first determine the differential of $\eta= \eta (\xi , x_2,\ldots ,x_n)$:
\begin{equation*}
d\eta=\frac{-2+\xi\left( \sum_{i=2}^{k}x_i^2-\sum_{i=k+1}^n x_i^2 \right)}{\xi^3}d\xi +\frac{ -2\sum_{i=2}^{k}x_i dx_i+2\sum_{i=k+1}^n x_i dx_i }{\xi}.
\end{equation*}
Hence, the pullback of $-\frac{1}{2}\partial^2 h$ to $\mathcal{H}_j$ which we denote by $g_j$, $j\in \{ 1,2\}$, is of the following form:
\begin{align*}
g_j=\ &\frac{3-\xi\left( \sum_{i=2}^{k}x_i^2-\sum_{i=k+1}^n x_i^2 \right)}{\xi^2} d\xi^2+2\left( \sum_{i=2}^{k}x_i dx_i - \sum_{i=k+1}^n x_i dx_i \right)d\xi\\
&+\xi \left( -\sum_{i=2}^k dx_i^2 + \sum_{i=k+1}^n dx_i^2 \right).
\end{align*}
For each $\frac{n+1}{2}\leq k\leq n$ there exists exactly one $\widetilde{k}$ with $1\leq\widetilde{k}\leq\frac{n+1}{2}$, such that $\mathcal{H}_1$ corresponding to $h=(x_1+x_{n+1})Q^{n+1}_k$ is isometric to $\mathcal{H}_2$ corresponding to $\widetilde{h}=(x_1+x_{n+1})Q^{n+1}_{\widetilde{k}}$, namely $\widetilde{k}=n-(k-1)$. In the coordinates $(\xi,\eta,x_2,\ldots,x_n)$ the corresponding isometry is given by $\xi\mapsto -\xi$, $x_\ell\mapsto x_{n-(\ell-2)}$ for $2\leq \ell\leq n$. Hence, we can reduce our analysis to $\mathcal{H}_1$, that is $\xi>0$, but need to increase the range for $k$ to $1\leq k\leq n$.

Returning to the study of $g_1$, we obtain
\begin{equation*}
g_1(\partial_{x_i},\partial_{x_j})=\left\{
\begin{array}{rl}
-\delta_i^j\xi, & 2\leq i,j\leq k, \\
\delta_i^j\xi, & k+1\leq i,j \leq n.
\end{array}\right.
\end{equation*}
For $g_1$ to be Riemannian, this implies that $k=1$. Hence, the only possibly hyperbolic polynomial is $h=\xi\left(\xi\eta-\sum_{k=2}^{n}x_i^2\right)$ and the corresponding metric $g_1$ reads
\begin{equation*}
g_1=\frac{3}{\xi^2}d\xi^2+\frac{1}{\xi}\sum_{i=2}^n\left( x_i d\xi -\xi dx_i\right)^2, 
\end{equation*}
which is indeed Riemannian at all points of $\mathcal{H}_1$. Hence, the only hyperbolic polynomial of the form $h=(x_1+x_{n+1})Q^{n+1}_k$, 
$\frac{n+1}{2}\leq k\leq n$, is given by
\begin{equation*}
h=(x_1+x_{n+1})\left( \sum_{i=1}^n x_i^2-x_{n+1}^2\right).
\end{equation*}
The corresponding $\mathcal{H}(h)=\{h=1,\ x_1+x_{n+1}<0\}$ has a single connected component. It is closed in $\mathbb{R}^{n+1}$, since $\{h=1\}\cap\{x_1+x_{n+1}=0\}=\emptyset$ implies that $\mathcal{H}(h)$ is also a connected component of $\{h=1\}$. This polynomial is the polynomial $b)$ of this proposition.
\end{proof}
\ep

\section{Curvature formulas for the q-map}
In this section, we introduce the supergravity r- and c-map and derive curvature formulas for their composition, the q-map. Note that compared to the last section, the dimension $n$ is shifted by one:
In this section, the projective special real manifold $\mathcal{H}$ is defined by a cubic polynomial $h$ in $n$ variables and has dimension $\mathrm{dim}\,\mathcal{H}=n-1$. The corresponding projective special K\"ahler manifold $\bar{M}$ in the image of the supergravity r-map has real dimension $2n$ and the quaternionic K\"ahler manifold $\bar{N}$ in the image of the q-map has real dimension $4m=4(n+1)$.
\subsection{Conical affine and projective special K\"ahler geometry}
\label{ConicalSect}
First, we recall the definitions of conical affine and projective special K\"ahler manifolds \cite{ACD,CM}:
\bd
A {\cmssl conical affine special K\"ahler manifold} $(M,g_M,J,\nabla,\xi)$ is a pseudo-K\"ahler manifold $(M,g_M,J)$ endowed with a flat torsionfree connection $\nabla$ and a vector field $\xi$ such that
\begin{enumerate}
\item[i)] $\nabla \omega_M=0$, where $\omega_M:=g_M(J\cdot,\cdot)$ is the K\"ahler form,
\item[ii)] $(\nabla_XJ)Y=(\nabla_YJ)X$ for all $X,Y\in \Gamma(TM)$,
\item[iii)] $\nabla\xi=D\xi=\mathrm{Id}$, where $D$ is the Levi-Civita connection,
\item[iv)] $g_M$ is positive definite on $\mathcal{D}=\mathrm{span}\{\xi,J\xi\}$ and negative definite on $\mathcal{D}^\perp$.
\end{enumerate}
\ed

Let $(M,J,g_M,\nabla,\xi)$ be a conical affine special K\"ahler manifold of complex dimension $n+1$. Then $\xi$ and $J\xi$ are commuting holomorphic vector fields that are homothetic and Killing respectively \cite{CM}. We assume that the holomorphic Killing vector field $J\xi$ induces a free $S^1$-action and that the holomorphic homothety $\xi$ induces a free $\mathbb{R}^{>0}$-action on $M$. Then $(M,g_M)$ is a metric cone over $(S,g_S)$, where $S:=\{p\in M|g_M(\xi(p),\xi(p))=1\}$, $g_S:=g_M|_S$; and $-g_S$ induces a Riemannian metric $g_{\bar{M}}$ on $\bar{M}:=S/S^1_{J\xi}$. $(\bar{M},-g_{\bar{M}})$ is obtained from $(M,J,g)$ via a 
K\"ahler reduction with respect to $J\xi$ and, hence, $g_{\bar{M}}$ is a K\"ahler metric (see e.g. \cite{CHM}). The corresponding K\"ahler form $\o_{\bar{M}}$ is
obtained from $\o_M$ by symplectic reduction. This determines
the complex structure $J_{\bar{M}}$. We will denote by $\pi$ the projection $M\to\bar{M}$. For future use let us mention that the metrics on $M$ and $\bar{M}$ are explicitly related by
\begin{equation}\label{eqn_gM_gbarM}
{g_M|}_{V\times V}=-g_M(\xi,\xi)\pi^* {g_{\bar{M}}|}_{V\times V},\quad V=(\ker d\pi)^\perp\subset TM.
\end{equation}

\bd
The K\"ahler manifold $(\bar{M},g_{\bar{M}},J_{\bar{M}})$ is called a {\cmssl projective special K\"ahler manifold}.
\ed

Locally, there exist so-called \emph{conical special holomorphic coordinates} $z=(z^I)=(z^0,\ldots ,z^n):U\stackrel{\sim}{\to} \tilde{U}\subset \mathbb{C}^{n+1}$ such that the geometric data on the domain $U\subset M$ is encoded in a holomorphic function $F: \tilde{U}\to \mathbb{C}$ that is homogeneous of degree 2 \cite{ACD,CM}. Namely, we have \cite{CM} \[g_M|_U=\sum_{I,J} N_{IJ}dz^Id\zb^J, \quad N_{IJ}(z,\zb):=2\mathrm{Im}\,F_{IJ}(z):=2\mathrm{Im}\,\frac{\partial^2F(z)}{\partial z^I \partial z^J} \quad (I,J=0,\ldots,n)\]
and $\xi|_U=\sum z^I\frac{\partial}{\partial z^I}+\zb^I\frac{\partial}{\partial\zb^I}$. The K\"ahler potential for $g_M|_U$ is given by $r^2|_U=g_M(\xi,\xi)|_U=\sum z^IN_{IJ} \zb^J$.

The $\mathbb{C}^\ast$-invariant functions $X^\mu:=\frac{z^\mu}{z^0}$, $\mu=1,\ldots, n$, define a local holomorphic coordinate system on $\bar{M}$. The K\"ahler potential for $g_{\bar{M}}$ is $\mathcal{K}:=-\log \sum_{I,J=0}^n X^IN_{IJ}(X)
\bar{X}^J$, where $X:=(X^0,\ldots, X^n)$ with $X^0:=1$. Note that for every function $f_U(z)$ on $U$, we define a function $f_{\bar{U}}(X)$ on the corresponding subset $\bar{U}\subset\bar{M}$ by $f_{\bar{U}}(X):=f_U(1,X^1,\ldots,X^n)$. In most cases, we will suppress the subscripts $_U$ and $_{\bar{U}}$ and use the same notation for corresponding functions on $U$ and $\bar{U}$.

\subsection{The supergravity c-map}
\label{sugraSect}
Let $(\bar{M},g_{\bar{M}})$ be a projective special K\"ahler manifold of complex dimension $n$ which is globally defined by a single holomorphic function $F$. The \emph{supergravity c-map} \cite{FS} associates with $(\bar{M},g_{\bar{M}})$ 
a quaternionic K\"ahler manifold $(\bar{N},g_{\bar{N}})$ of dimension $4n+4$. 
Following the conventions of \cite{CHM}, we have $\bar{N}=\bar{M}\times \mathbb{R}^{>0}\times \mathbb{R}^{2n+3}$ and
\begin{eqnarray*} \label{e:2} g_{\bar{N}} &=& g_{\bar{M}} + g_G,\\\nonumber 
g_G&=& \frac{1}{4\rho^2}d\rho^2 + \frac{1}{4\rho^2}(d\tilde{\phi}
+ \sum (\zeta^Id\tilde{\zeta}_I-\tilde{\zeta}_Id\zeta^I) )^2 
+\frac{1}{2\rho}\sum \mathcal{I}_{IJ}(m) d\zeta^Id\zeta^J\\ 
&&+ \frac{1}{2\rho}\sum \mathcal{I}^{IJ}(m)(d\tilde{\zeta}_I + 
\mathcal{R}_{IK}(m)d\zeta^K)(d\tilde{\zeta}_J + 
\mathcal{R}_{JL}(m)d\zeta^L),
\end{eqnarray*}
where $(\rho,\tilde{\phi},\tilde{\z}_I, \z^I)$, $I=0,1,\ldots, n$, 
are standard coordinates on $\bR^{>0}\times \bR^{2n+3}$. 
The real-valued matrices $\mathcal{I}(m):=(\mathcal{I}_{IJ}(m))$ and $\mathcal{R}(m):=(\mathcal{R}_{IJ}(m))$
depend only on $m\in \bar{M}$ and $\mathcal{I}(m)$ is invertible with
the inverse $\mathcal{I}^{-1}(m)=:(\mathcal{I}^{IJ}(m))$. More precisely,
 \begin{equation*} 
\label{FRI}
{\cal N}_{IJ} := 
\mathcal{R}_{IJ} + i\mathcal{I}_{IJ} := 
\bar{F}_{IJ} + 
i \frac{\sum_K N_{IK}z^K\sum_L N_{JL}z^L}{\sum_{IJ}N_{IJ}z^Iz^J} ,\quad 
N_{IJ} := 2 \,\mathrm{Im} \,F_{IJ} ,\end{equation*}
where $F$ is the holomorphic prepotential with respect
to some system of special holomorphic coordinates $z^I$ on the 
underlying conical special K\"ahler manifold $M\ra \bar{M}$. 
Notice that the expressions are homogeneous of degree zero and, hence, 
well defined functions on $\bar{M}$. It is shown in \cite[Cor.\ 5]{CHM} 
that the matrix $\mathcal{I}(m)$ is positive definite
and hence invertible and that the metric $g_{\bar{N}}$ does not
depend on the choice of special coordinates \cite[Thm.\ 9]{CHM}. 
It is also shown that $(\bar{N},g_{\bar{N}})$ is complete if 
and only if $(\bar{M},g_{\bar{M}})$ is complete \cite[Thm.\ 5]{CHM}.

Using $(p_a)_{a=1,\ldots, 2n+2}:=(\zt_I,\z^J)_{IJ=0,\ldots,n}$ and $(\hat{H}^{ab}):=\begin{pmatrix}\mathcal{I}^{-1} & \mathcal{I}^{-1}\mathcal{R} \\ \mathcal{R}\mathcal{I}^{-1} & \mathcal{I}+\mathcal{R}\mathcal{I}^{-1}\mathcal{R}\end{pmatrix}$, we can combine the last two terms of $g_G$ into $\frac{1}{2\rho}\sum dp_a \hat{H}^{ab} dp_b$, i.e. the quaternionic K\"ahler metric is given by
\begin{equation}\label{FSmetric} g_{FS}:=g_{\bar{N}}=g_{\bar{M}}+\frac{1}{4\rho^2}d\rho^2 + \frac{1}{4\rho^2}(d\tilde{\phi}
+ \sum (\zeta^Id\tilde{\zeta}_I-\tilde{\zeta}_Id\zeta^I) )^2+ \frac{1}{2\rho}\sum dp_a \hat{H}^{ab} dp_b.\end{equation}

\subsection{The supergravity r-map}\label{secSugraRMap}
Let $(\mathcal{H}:=\{x\in U\mid h(x)=1\}, g_\mathcal{H}:=-\partial^2h\big|_{\mathcal{H}})$ be a projective special real manifold defined by a real homogeneous cubic polynomial $h$ and an $\mathbb{R}^{>0}$-invariant domain $U\subset \mathbb{R}^n\backslash\{0\}$. Let $\bar{M}:=\mathbb{R}^n+iU\subset \mathbb{C}^n$ be endowed with the standard complex structure $J_{\bar{M}}$ induced from $\mathbb{C}^n$ and with holomorphic coordinates 
$(X^\mu=y^\mu+ix^\mu)_{\mu=1,\ldots,n}\in\mathbb{R}^n+iU$. We define a K\"ahler metric
\begin{equation*}
g_{\bar{M}}=\sum_{\mu,\nu=1}^n\frac{\partial^2\mathcal{K}}{\partial X^\mu\partial\bar{X}^\nu} dX^\mu d\bar{X}^\nu
\end{equation*}
on $\bar{M}$ with K\"ahler potential
\begin{equation*}
\mathcal{K}(X,\bar{X}):=-\log 8h(x),
\end{equation*}
where $x=(\opn{Im} X^1,\ldots,\opn{Im}X^n)\in U$.

\bd
The correspondence $(\mathcal{H},g_{\mathcal{H}})\mapsto (\bar{M},g_{\bar{M}},J_{\bar{M}})$ is called the {\bf supergravity r-map}.
\ed

\br\label{remHolPrepotrMap}{\rmfamily\normalfont
Note that any manifold $(\bar{M},\,g_{\bar{M}},\,J_{\bar{M}})$ in the image of the supergravity r-map is a projective special K\"ahler manifold (see Section \ref{ConicalSect}). The corresponding conical affine special K\"ahler manifold is the trivial $\mathbb{C}^\ast$-bundle \[M:=\{z=z^0\cdot(1,\,X)\in\mathbb{C}^{n+1}\mid z^0\in\mathbb{C}^\ast,~X\in\bar{M}=\mathbb{R}^n+iU\}\to\bar{M}\] endowed with the standard complex structure $J$ and the metric $g_M$ defined by the holomorphic function
\[F:M\to \mathbb{C},\quad F(z^0,\,\ldots,\,z^n)=\frac{h(z^1,\,\cdots,\,z^n)}{z^0}.\]
Note that in general, the flat connection\footnote{$\nabla$ is defined by $\mathbf{x}^I=\opn{Re} z^I$ and $\mathbf{y}_I=\opn{Re} F_I(z)$ being flat for $I=0,\,\ldots,\,n$ (see \cite{ACD}).} $\nabla$ on $M$ is not the standard one induced from $\mathbb{C}^{n+1}\approx \mathbb{R}^{2n+2}$. The homothetic vector field $\xi$ is given by $\xi=\sum_{I=0}^n(z^I\frac{\partial}{\partial z^I}+\bar{z}^I\frac{\partial}{\partial \bar{z}^I})$. To check that $g_{\bar{M}}$ is the corresponding projective special K\"ahler metric, one uses the fact that
\begin{equation*}
8|z^0|^2h(x)=\sum_{I,\,J=0}^nz^IN_{IJ}(z,\,\bar{z})\bar{z}^J,
\end{equation*}
where as above, $x=(\opn{Im} X^1,\,\ldots,\,\opn{Im} X^n)=(\opn{Im} \frac{z^1}{z^0},\,\ldots,\,\opn{Im} \frac{z^n}{z^0})\in U$ (see \cite{CHM}).
}\er

\subsection{Curvature formulas for the supergravity r-map}
Under the assumptions of Section \ref{secSugraRMap}, let $(e^a_\mu)_{a,\,\mu=1,\,\ldots,\,n}$ be a real $n\times n$ matrix-valued function on some open subset in $\bar{M}$ such that $\sum_{a=1}^n e^a_\mu \bar{e}^a_\nu=\sum_{a=1}^n e^a_\mu e^a_\nu=\mathcal{K}_{\mu\bar{\nu}}$, where
\be\label{eqSugraRMapMetric} \mathcal{K}_{\mu\bar{\nu}}=-\frac{\partial^2\log\,h(x)}{\partial X^\mu\partial{\bar{X}}^\nu}=-\frac{h_{\mu\nu}(x)}{4h(x)}+\frac{h_\mu(x)h_\nu(x)}{4h^2(x)}.\ee
Here, subscripts of the cubic polynomial $h$ denote derivatives with respect to the standard coordinates on $U$, e.g.\ $h_\mu(x)=\frac{\partial h(x)}{\partial x^\mu}$.
The holomorphic one-forms
\be\label{sigmas} \sigma^a:=\sum_{\mu=1}^n e^a_\mu dX^\mu\ee
constitute a \emph{unitary coframe} $(\sigma^a)_{a=1,\,\ldots,\,n}$, i.e.\ the metric can locally be written as
\begin{equation*}
g_{\bar{M}}=\sum_{a=1}^n\sigma^a\bar{\sigma}^a=\frac{1}{2}\sum_{a=1}^n(\sigma^a\otimes\bar{\sigma}^a+\bar{\sigma}^a\otimes\sigma^a).
\end{equation*}
Let $(\sigma_a:=\sum_{\mu=1}^n e^\mu_a\frac{\partial}{\partial X^\mu})_{a=1,\,\ldots,\,n}$ denote the corresponding local frame in $T^{1,\,0}\bar{M}$ dual to $(\sigma^a)_{a=1,\,\ldots,\,n}$, i.e.\ $(e^\mu_a)=(e^a_\mu)^{-1}$. Then $\sigma^a=2g_{\bar{M}}(\bar{\sigma}_a,\,\cdot)$ and $\sigma^a(\sigma_b)=\bar{\sigma}^a(\bar{\sigma}_b)=\delta^a_{~b}$, \linebreak$\sigma^a(\bar{\sigma}_b)=\bar{\sigma}^a(\sigma_b)=0$. Note that $g_{\bar{M}}(\sigma_a,\bar{\sigma}_b)=\frac{1}{2}\delta_{ab}$ which implies $\langle \sigma_a,\sigma_b\rangle=\delta_{ab}$ for the corresponding sesqui-linear form $\langle\cdot,\cdot\rangle=g_{\bar{M}}-i\omega_{\bar{M}}$ on $(T\bar{M},J)$. This explains why we call the frame $(\sigma_a)$ and the dual coframe $(\sigma^a)$ unitary.

\noindent
Note that the inverse of the matrix-valued function $(\mathcal{K}_{\mu\bar{\nu}})_{\mu,\nu=1,\ldots,n}$ (see equation \eqref{eqSugraRMapMetric}) is given by
\be\label{eqn_Kinverse} \mathcal{K}^{\bar{\nu}\rho}=-4h(x)h^{\nu\rho}(x)+2x^\nu x^\rho,\ee
where $(h^{\mu\nu})_{\mu,\nu=1,\ldots,n}=(h_{\mu\nu})_{\mu,\nu=1,\ldots,n}^{-1}$. We will usually write $\mathcal{K}^{\rho\bar{\nu}}$ instead of $\mathcal{K}^{\bar{\nu}\rho}$.

\noindent
Note that in this section, $\nabla$ denotes the Levi-Civita connection of the projective special K\"ahler metric $g_{\bar{M}}$.
The expressions for the Christoffel symbols
\begin{align*}\nonumber \Gamma^\rho_{\sigma\mu}&:=dX^\rho(\nabla_{\partial_{X^\sigma}}\partial_{X^\mu})=\sum_{\kappa=1}^n\mathcal{K}^{\rho\bar{\kappa}}\partial_{X^\sigma} \mathcal{K}_{\mu\bar{\kappa}}\\
& ~=-\frac{i}{2h}\left(h\sum_{\kappa=1}^n h^{\rho\kappa}h_{\kappa\mu\sigma}-h_\sigma \delta^\rho_\mu-h_\mu\delta^\rho_\sigma +\frac{1}{2}x^\rho h_{\mu\sigma}\right)
\end{align*}
and the coefficients
\begin{align}
\nonumber R^\rho_{~\sigma\mu\bar{\nu}}&:=dX^\rho\left( R(\partial_{X^\mu},\,\partial_{{\bar{X}}^{\nu}}) \partial_{X^\sigma}\right)=-\partial_{\bar{X}^\nu}\Gamma^\rho_{\sigma\mu}=-\frac{i}{2}\partial_{x^\nu}\Gamma^\rho_{\sigma\mu}\\
&~=-\frac{1}{4h^2}\Bigg[ \frac{1}{2}x^\rho(h h_{\mu\sigma\nu}-h_{\mu\sigma}h_\nu)+h_\mu h_\nu\delta_\sigma^\rho+h_\sigma h_\nu\delta_\mu^\rho \nonumber \\
&~~\qquad\qquad -h\left( h_{\sigma\nu}\delta^\rho_\mu+h_{\mu\nu}\delta_\sigma^\rho-\frac{1}{2}h_{\mu\sigma}\delta_\nu^\rho\right)-h^2\sum_{\alpha,\beta,\gamma=1}^n h^{\rho \alpha}h_{\nu\alpha\beta}h^{\beta\gamma}h_{\gamma\mu\sigma} \Bigg]\nonumber\\&=-\delta^\rho_\sigma\mathcal{K}_{\mu\bar{\nu}}-\delta^\rho_\mu\mathcal{K}_{\sigma\bar{\nu}}+e^{2\mathcal{K}}\sum_{\alpha,\beta,\gamma=1}^n\mathcal{K}^{\rho\bar{\alpha}}h_{\alpha\nu\beta}\mathcal{K}^{\beta\bar{\gamma}}h_{\gamma\mu\sigma} \label{curvatureTensor}
\end{align}
of the Riemann curvature tensor
\begin{equation*}\quad\quad R(X,Y)Z:=\nabla_X\nabla_Y Z-\nabla_Y\nabla_XZ-\nabla_{[X,Y]}Z \quad\quad(X,Y,Z\in \mathfrak{X}(\bar{M}))\end{equation*}
have been calculated for instance in \cite[Theorem 3]{CDL}.

\noindent We denote the coefficients of the local Levi-Civita connection one-form associated to the unitary local coframe $(\sigma^a)_{a=1,\ldots,n}$ by $\omega^a_{~b}$, i.e.\ $\nabla_\cdot \sigma^a=-\sum_{b=1}^n \omega^a_{~b}(\cdot)\sigma^b$. Compatibility with the metric and torsion-freeness translate into the conditions that the complex one-form valued matrix $(\omega^a_{~b})_{a,b=1,\ldots,n}$ is anti-Hermitian and satisfies $d\sigma^a+\sum_{b=1}^n\omega^a_{~b}\wedge\sigma^b=0$ for $a=1,\ldots,n$. These are fulfilled by the following general formula that holds for all K\"ahler manifolds:
\footnote{Note that for arbitrary K\"ahler manifolds, the functions $e^a_\mu$ cannot in general be chosen to be real.}
\begin{align}
\label{eqConnectionOneFormKahlerGeneral} \omega^a_{~b}&=\sum_{\mu=1}^n(e^a_\mu \bar{\partial}e^\mu_b-\bar{e}^b_\mu\partial\bar{e}^\mu_a).
\end{align}
This formula is found by observing that the $(0,1)$-component of $\omega^a_{~b}$ is uniquely determined by solving the $(1,1)$-projection of the equation $d\sigma^a+\sum_{b=1}^n\omega^a_{~b}\wedge\sigma^b=0$ and using the skew-Hermiticity to compute the $(1,0)$-component of $\omega^a_{~b}$. By the existence of the Levi-Civita connection the $(2,0)$-projection of the equation $d\sigma^a+\sum_{b=1}^n\omega^a_{~b}\wedge\sigma^b=0$ is then automatically satisfied.
In terms of the local connection one-form, the curvature tensor of a K\"ahler manifold is given by
\begin{equation*}
R(X,\,Y)\sigma_c=\sum_{d=1}^n(d\omega^d_{~c}+\sum_{c'=1}^n\omega^d_{~c'}\wedge\omega^{c'}_{~c})(X,\,Y)\sigma_d=:\sum_{d=1}^n\tilde{R}^d_{~c}(X,\,Y)\sigma_d.
\end{equation*}
Using equation \eqref{curvatureTensor} and $\mathcal{K}^{\mu\bar{\nu}}=\sum_{c=1}^n e^\mu_c\bar{e}^\nu_c$, one gets the following proposition (see \cite[Prop.\ 7.2.1]{D}):
\bp
In terms of the unitary local coframe $(\sigma^a)_{a=1,\,\ldots,\,n}$, the Riemann curvature tensor of a projective special K\"ahler manifold in the image of the supergravity r-map reads
\begin{equation*}
\tilde{R}^a_{~b}=-\delta^a_b\sum_{c=1} ^n\sigma^c\wedge\bar{\sigma}^c-\sigma^a\wedge\bar{\sigma}^b+e^{2\mathcal{K}}\sum_{c,e,d=1}^n\tilde{h}_{adc}\tilde{h}_{ceb}\sigma^e\wedge\bar{\sigma}^d,
\end{equation*}
where $\tilde{h}_{abc}:=\sum_{\mu,\nu,\sigma=1}^n e^\mu_a e^\nu_b e^\sigma_c h_{\mu\nu\sigma}$ for $a,b,c=1,\ldots, n$.
\ep
\begin{proof}
Using $\mathcal{K}_{\mu\bar{\nu}}=e^c_\mu\bar{e}^c_\nu$, $\mathcal{K}^{\mu\bar{\nu}}=e_c^\mu\bar{e}_c^\nu$ and the fact that $(e^a_\mu)_{a,\,\mu=1,\,\ldots,\,n}=(e^\nu_b)_{\nu,\,b=1,\,\ldots,\,n}^{-1}$, we find
\begin{align*}\tilde{R}^a_{~b}(\sigma_e,\,\bar{\sigma}_d)&=
\sigma^a(R(\sigma_e,\,\bar{\sigma}_d)\sigma_b)\\&=e^a_\rho \,dX^\rho\!\left( R(\partial_{X^\mu},\,\partial_{{\bar{X}}^{\nu}}) \partial_{X^\sigma}\right)e^\mu_e\bar{e}^\nu_de^\sigma_b\\
&\!\!\!\stackrel{\eqref{curvatureTensor}}{=}e^a_\rho \left(-\delta^\rho_\sigma\mathcal{K}_{\mu\bar{\nu}}-\delta^\rho_\mu\mathcal{K}_{\sigma\bar{\nu}}+e^{2\mathcal{K}}\mathcal{K}^{\rho\bar{\alpha}}h_{\alpha\nu\beta}\mathcal{K}^{\beta\bar{\gamma}}h_{\gamma\mu\sigma}\right)e^\mu_e\bar{e}^\nu_de^\sigma_b\\
&=-\delta^a_b\delta_{de}-\delta^a_e\delta_{bd}+e^{2\mathcal{K}}\tilde{h}_{adc}\tilde{h}_{ceb}.
\end{align*}
\end{proof}

\subsection{Levi-Civita connection for quaternionic K\"ahler manifolds obtained by the q-map}
In this and the following section, we will introduce the \emph{quaternionic vielbein formalism}, which was used in \cite{FS} to determine the Levi-Civita connection and the Riemann curvature tensor of manifolds in the image of the supergravity c-map. The formulas in this formalism arise from well-known formulas in the differential geometry literature expressed in terms of local frames in the complex vector bundles $E$ and $H$ whose tensor product is identified with the complexified tangent bundle of a quaternionic K\"ahler manifold in Salamon's $E$-$H$ formalism \cite{S} (see e.g.\ \cite[Ch.\ 7]{D} for detailed explanations of the relation between the formulas used in the physics, respectively mathematics literature). The \emph{q-map} is the composition of the supergravity r- and c-map. It assigns a quaternionic K\"ahler manifold of dimension $4m=4(n+1)$ to any projective special real manifold of dimension $n-1$. We apply the quaternionic vielbein formalism to quaternionic K\"ahler manifolds in the image of the q-map and derive formulas for the Levi-Civita connection and the Riemann curvature tensor of these manifolds, expressed in terms of the cubic polynomial $h$, which defines the initial projective special real manifold. Up to changing conventions and fixing inaccuracies, these results can also be obtained by restricting the formulas in \cite{FS} for the c-map to the case of the q-map. The Riemann curvature tensor of a quaternionic K\"ahler manifold is determined by its trace-free part, the \emph{quaternionic Weyl tensor}. The latter can be expressed in terms of a certain symmetric quartic tensor field $\Omega\in \Gamma(S^4E^\ast)$ in the complex vector bundle $E$.
In addition to the above-mentioned results, we derive a formula expressing this quartic tensor field in terms of the cubic polynomial $h$ for manifolds in the image of the q-map. This result is used in Subsection \ref{subsecNormCurv} to give a general formula for the squared pointwise norm of the Riemann curvature tensor of any quaternionic K\"ahler manifold in the image of the q-map.

We will restrict ourselves to manifolds in the image of the \emph{q-map}, which is the composition of the supergravity r- and c-map, i.e.\ we consider the Ferrara-Sabharwal metric \eqref{FSmetric} defined on $\bar{N}=\bar{M}\times \mathbb{R}^{>0}\times \mathbb{R}^{2n+3}$ for a projective special K\"ahler manifold $(\bar{M}=\mathbb{R}^n+iU,g_{\bar{M}},J_{\bar{M}})$ in the image of the supergravity r-map, which is defined by a real homogeneous cubic polynomial $h$. On $\bar{N}$, we define the following complex-valued one-forms:
\begin{align}\label{eqQuaternionicCoframe}
\beta^0&:=i e^{\mathcal{K}/2}\frac{1}{\sqrt{\rho}}\sum_{I=0}^n X^IA_I,&
\beta^a&:=\sum_{I=0}^n P^a_IdX^I=\sigma^a,\\\nonumber
\alpha^0&:=-\frac{1}{2\rho}\left(d\rho-i(d\tilde{\phi}+\sum_{I=0}^n(\zeta^Id\zt_I-\zt_Id\z^I))\right),&
\alpha^a&:=\frac{i}{\sqrt{\rho}}e^{-\mathcal{K}/2}\sum_{I,J=0}^n\overline{P}^a_IN^{IJ}A_J
\end{align}
for $a=1,\ldots,n$, where $X^0=1$, $(N^{IJ})$ is the inverse of the matrix $(N_{IJ})$,
$P^a_I$ are the components of the complex $n\times(n+1)$ matrix-valued function \be\label{PI_eqn}(P^a_I)_{a=1,\,\ldots,\,n,~I=0,\,\ldots,\,n}=(P^a_0,\,P^a_\mu)_{a,\,\mu=1,\,\ldots,\,n}:=\left(-\sum_{\nu=1}^ne^a_\nu X^\nu,\,e^a_\mu\right)_{a,\,\mu=1,\,\ldots,\,n},\ee
and $A_I=d\zt_I+\sum_{J=0}^n F_{IJ}(X)d\z^J$ for $I=0,\ldots,n$. Note that the matrix $(P^a_I)$ represents the linear map $d\pi_p:T_pM\to T_{\pi(p)}\bar{M}$ for $p=(X^0,X^1,\ldots,X^n)^t$ with $X^0=1$ in the coordinate basis $\left(\frac{\partial}{\partial X^I}\right)$ of $T_pM$ and the unitary basis $(\sigma_a)$ of $T_{\pi(p)}\bar{M}$. In terms of these one-forms, the Ferrara-Sabharwal metric reads (see e.g.\ \cite[Lemma 7.3.1]{D})
\begin{equation*}
g_{FS}=\sum_{A=0}^n(\beta^A\bar{\beta}^A+\alpha^A\bar{\alpha}^A).
\end{equation*}
The equations
\begin{equation}\label{eqn_J1J2J3}
J_1^\ast \alpha^A=i\alpha^A,\quad J_1^\ast \beta^A=i\beta^A,\quad J_2^\ast\alpha^A=\bar{\beta}^A,\quad 
\end{equation}
for $A=0,\ldots,n$ and $J_1J_2=J_3$ define an almost hyper-complex structure $(J_1,J_2,J_3)$ on $\bar{N}$. $J_1$, $J_2$ and $J_3$ span a quaternionic structure $Q$ on $\bar{N}$ that is compatible (skew-symmetric and parallel) with the quaternionic K\"ahler metric $g_{FS}$. In fact, the quaternionic K\"ahler property was proven in \cite{FS} by computing the Levi-Civita connection. This calculation is reviewed below, see equations \re{eqn_LCqk}--\re{eqQKCOnnE} and Proposition \ref{cor_SPn}, and amounts to showing that $Q$ is parallel. Alternatively, it follows from the fact that the Ferrara-Sabharwal metric can be obtained by the geometric construction described in \cite{ACDM} involving the HK/QK-correspondence. Note that $J_1$ defines an integrable\footnote{This can either be shown by direct calculation (see \cite{CLST}) or deduced from the fact that all quaternionic K\"ahler manifolds obtained from the HK/QK correspondence admit a globally defined compatible integrable complex structure (see \cite[Rem.\ 5.5.5]{D}).} complex structure on $\bar{N}$.

We will now prepare for the calculation of the Levi-Civita connection of the Ferrara-Sabharwal metric. Direct calculation gives the following expressions for the exterior derivatives of the above one-forms (see \cite[Prop.\ 7.3.3]{D}):
\bp\label{prop_differentials}
\begin{align*}
d\beta^0&=\frac{1}{2}\left(\alpha^0+\bar{\alpha}^0-id^c\mathcal{K}\right)\wedge \beta^0+\sum_{b=1}^n \alpha^b\wedge \beta^b,\\
d\beta^a&=-\sum_{b=1}^n\omega^a_{~b}\wedge \beta^b,\\
d\alpha^0&=-\alpha^0\wedge\bar{\alpha}^0+\beta^0\wedge \bar{\beta}^0-\sum_{b=1}^n\alpha^b\wedge \bar{\alpha}^b,
\end{align*}
\begin{align*}
d\alpha^a&=\frac{1}{2}(\alpha^0+\bar{\alpha}^0-id^c\mathcal{K})\wedge \alpha^a+\beta^0\wedge \bar{\beta}^a-\sum_{b=1}^n\overline{\omega^a_{~b}}\wedge \alpha^b-ie^{\mathcal{K}}\sum_{b,c=1}^n\tilde{h}_{abc}\bar{\alpha}^b\wedge \beta^c,
\end{align*}
where $d^c=i(\bar{\partial}-\partial)$, $\tilde{h}_{abc}=\sum_{\mu,\nu,\sigma=1}^n e^\mu_a e^\nu_b e^\sigma_c h_{\mu\nu\sigma}$ for $a,b,c=1,\ldots, n$, and $(\omega^a_{~b})_{a,b=1,\ldots,n}$ is the (pullback to $\bar{N}$ of the) local connection one-form of the Levi-Civita connection on $\bar{M}$ with respect to the local unitary coframe $(\sigma^a)_{a=1,\ldots,n}$ on $\bar{M}$. \label{propositionDerivativeCoframeqMap}
\ep

\noindent
To calculate the exterior derivatives in Proposition \ref{prop_differentials} we have used the following explicit formula for the local Levi-Civita connection one-form of a projective special K\"ahler manifold:
\bp\label{propRMapConnectionOneForm}
The local connection one-form for the Levi-Civita connection with respect to the unitary coframe $(\sigma^a)_{a=1,\,\ldots,\,n}$ can be written as
\begin{align*}
\omega^a_{~b}&=e^{-\mathcal{K}}\big((\bar{\partial}P^a_I)N^{IJ}\bar{P}^b_J-P^a_I N^{IJ}(\partial\bar{P}^b_J)\big)
\\&=\delta^a_b\partial\mathcal{K}+e^{-\mathcal{K}}d(P^a_IN^{IJ})\bar{P}^b_J+ie^{-\mathcal{K}}P^a_IN^{IK}\,\overline{dF_{KL}(X)}\,N^{LJ}\bar{P}^b_J,
\end{align*}
where the $P^a_I$ are defined in equation \re{PI_eqn}.
\ep
\begin{proof}
Note that ${-e^{-\mathcal{K}}g_{\bar{M}}(d\pi\cdot,d\pi\cdot)|}_{V\times V}={g_M|}_{V\times V}$, see \re{eqn_gM_gbarM} and $e^{-\mathcal{K}}={g_M(\xi,\xi)|}_{\{X^0=1\}}$. Dualizing yields the equation
\begin{equation*}
	-e^{-\mathcal{K}}g^{-1}_{\bar{M}}=g^{-1}_M(d\pi^*\cdot,d\pi^*\cdot),
\end{equation*}
which in components reads
\begin{align}
\label{eqPNPEqDelta}
-e^\mathcal{K}\delta^{ab}&=
\sum_{I,\,J=0}^nP^a_IN^{IJ}\bar{P}^b_J& &(a,\,b=1,\,\ldots,\,n).\end{align}
Multiplication of equation \eqref{eqPNPEqDelta} by $-e^{-\mathcal{K}}e_a^\mu$ gives
\[ e^\mu_b=e^{-\mathcal{K}}(X^\mu N^{0J}-N^{\mu J})\bar{P}^b_J.\]
This equation shows that
\[-e^\mu_b\bar{\partial}e^a_\mu=e^{-\mathcal{K}}((\bar{\partial}e^a_\mu)N^{\mu J}-X^\mu(\bar{\partial}e^a_\mu)N^{0J})\bar{P}^b_J=e^{-\mathcal{K}}(\bar{\partial}P^a_I)N^{IJ}\bar{P}^b_J.\]
Using the above equation one then finds
\[\omega^a_{~b}\stackrel{\eqref{eqConnectionOneFormKahlerGeneral}}{=}-e^\mu_b\bar{\partial} e^a_\mu+\bar{e}^\mu_a\partial \bar{e}^b_\mu=e^{-\mathcal{K}}\big((\bar{\partial}P^a_I)N^{IJ}\bar{P}^b_J-P^a_I N^{IJ}(\partial\bar{P}^b_J)\big).\]
Adding $0\stackrel{\eqref{eqPNPEqDelta}}{=}\delta^a_b\partial\mathcal{K}+e^{-\mathcal{K}}\partial(P^a_IN^{IJ}\bar{P}^b_J)$ to the above equation gives
\begin{align*}
\omega^a_{~b}&=\delta^a_b\partial\mathcal{K}+e^{-\mathcal{K}}\big(d(P^a_IN^{IJ})\bar{P}^b_J-P^a_I(\bar{\partial}N^{IJ})\bar{P}^b_J\big)\\
&=\delta^a_b\partial\mathcal{K}+e^{-\mathcal{K}}d(P^a_IN^{IJ})\bar{P}^b_J+ie^{-\mathcal{K}}P^a_IN^{IK}\,\overline{dF_{KL}(X)}\,N^{LJ}\bar{P}^b_J.
\end{align*}
\end{proof}

The components $\bar{\theta}_\alpha$ of the local $\mathrm{Sp}(1)$-connection one-form of a quaternionic K\"ahler manifold (with Levi-Civita connection $\nabla$) with respect to a local oriented orthonormal frame $(J_1,J_2,J_3)$ in the quaternionic structure are defined by
\begin{equation*}
\nabla_\cdot J_\alpha=2(\bar{\theta}_\beta(\cdot)J_\gamma-\bar{\theta}_\gamma(\cdot)J_\beta)
\end{equation*}
for any cyclic permutation $(\alpha,\beta,\gamma)$ of $(1,2,3)$.
The one-forms $\bar{\theta}_\alpha$ are related to the fundamental two-forms $\omega_\alpha=g(J_\alpha\cdot,\cdot)$ by the following well known structure equations for quaternionic K\"ahler manifolds
\be\label{propChQKAlekseevskyCurvature} \frac{\nu}{2}\omega_\alpha=d\bar{\theta}_\alpha-2\bar{\theta}_\beta\wedge\bar{\theta}_\gamma,\ee
where $\nu:=\frac{scal}{4m(m+2)}$ ($\mathrm{dim}_\mathbb{R}\bar{N}=4m=4(n+1)$) is the \emph{reduced scalar curvature}. For manifolds in the image of the supergravity c-map, we have $\nu=-2$ and 
\begin{align}\label{eqn_theta123}
\bar{\theta}_1&=-\frac{1}{4\rho}\big(d\tilde{\phi}+\rho \,d^c\mathcal{K}-\sum_{I=0}^n(\tilde{\zeta}_Id\zeta^I-\zeta^I d\tilde{\zeta}_I)\big)=-\frac{1}{2}\mathrm{Im}\,\alpha^0-\frac{1}{4}d^c\mathcal{K},\nonumber\\
\bar{\theta}_2+i\bar{\theta}_3&=i\frac{1}{\sqrt{\rho}}e^{\mathcal{K}/2} \sum_{I=0}^n X^I A_I=\beta^0.
\end{align}
These formulas follow from the general expression for the $\mathrm{Sp}(1)$-connection of a quaternionic K\"ahler manifold obtained from the HK/QK-correspondence, see \cite[Thm.\ 4.1.2]{D}, after specialization to the case of the supergravity c-map, see \cite[Rem.\ 5.5.3 and 5.5.4]{D} and \cite{ACDM}.

We combine the one-forms defined in equation \eqref{eqQuaternionicCoframe} into the following \emph{quaternionic vielbein}, which is a $(4n+4)\times(4n+4)$ matrix of complex-valued one-forms:
\begin{equation*}
(f^{\alpha\Gamma})_{\alpha=1,2;\Gamma=1,\ldots,2n+2}=\begin{pmatrix}
f^{1A} & f^{1\tilde{A}}\\f^{2A} & f^{2\tilde{A}}
\end{pmatrix}_{A=0,\ldots,n}:=\begin{pmatrix}\beta^A & \alpha^{A}\\-\bar{\alpha}^A& \bar{\beta}^{A}\end{pmatrix}_{A=0,\ldots,n}.
\end{equation*}
Let $\beta_A$, $\alpha_A$ be complex-valued vector fields on $\bar{N}$ such that $\beta^A=2\overline{g(\beta_A,\cdot)}$ and\linebreak $\alpha^A=2\overline{g(\alpha_A,\cdot)}$ for $A=0,\ldots,n$. These vector-fields are combined into the following local frame in $T^\mathbb{C}\bar{N}$, which is dual to $(f_{\alpha\Gamma})$:
\be\label{eqQuaternionicVielbeinFrame} (f_{\alpha\Gamma})_{\alpha=1,2;\Gamma=1,\ldots,2n+2}=\begin{pmatrix}
f_{1A}&f_{1\tilde{A}}\\f_{2A}&f_{2\tilde{A}}\end{pmatrix}_{A=0,\ldots,n}:=\begin{pmatrix}
\beta_A&\alpha_A\\-\bar{\alpha}_A&\bar{\beta}_A
\end{pmatrix}_{A=0,\,\ldots,\,n}. \ee
Recall that the skew-symmetric almost complex structures $J_1,J_2,J_3$ spanning the quaternionic structure $Q$ are of standard form, see \re{eqn_J1J2J3}, in the coframe $(f^{\alpha\Gamma})$. Note that in our case, namely for manifolds in the image of the q-map, the frame $(f_{\alpha\Gamma})$ is globally defined and thus establishes a global isomorphism $T^\mathbb{C}\bar{N}\cong H\otimes E$, where $H$ and $E$ are trivial complex vector bundles. More specifically, $(f_{\alpha\Gamma})$ corresponds to a tensor product of the form $(f_{\alpha\Gamma})=(h_\alpha\otimes E_\Gamma)$, where $(h_\alpha)$ is a frame of $H$ and $(E_\Gamma)$ is a frame of $E$.
To prove that $Q$ is parallel and therefore that $g_{FS}$ is quaternionic K\"ahler, it is sufficient to check that the Levi-Civita connection with respect to frame $(f_{\alpha\Gamma})$ has the following form:
\be\label{eqn_LCqk} f^{\alpha \Gamma}(\nabla_Xf_{\beta\Delta})=p^\alpha_{~\beta}(X)\delta^\Gamma_{~\Delta}+\delta^\alpha_{~\beta}\Theta^\Gamma_{~\Delta}(X)\ee
for $\alpha,\beta=1,2$ and $\Gamma,\Delta=1,\ldots,2n+2$,
where $p=(p^\alpha_{~\beta})$ is a one-form with values in $\mathfrak{sp}(1)=\mathfrak{su}(2)$, i.e. $p^\dagger:=\bar{p}^{t}=-p$, and $\Theta=(\Theta^{\Gamma}_{~\Delta})$ is a one-form with values in $\mathfrak{sp}(n+1)\subset\mathfrak{su}(2n+2)$. The latter means that
\be \Theta=\begin{pmatrix}q&t\\-\bar{t}&\bar{q}\end{pmatrix},\label{eqQKCOnnE}\ee
where $q$, $t$ are complex 1-form-valued $(n+1)\times (n+1)$ matrices that are anti-Hermitian, respectively symmetric ($q^\dagger=\bar{q}^t=-q$, $t^t=t$).
\bp\label{cor_SPn}
The $\mathrm{Sp}(1)$-part of the Levi-Civita connection of a quaternionic K\"ahler manifold in the image of the q-map is given by
\begin{equation*}
p=\begin{pmatrix}-i\bar{\theta}_1&-\bar{\theta}_2-i\bar{\theta}_3\\\bar{\theta}_2-i\bar{\theta}_3&i\bar{\theta}_1\end{pmatrix},
\end{equation*}
see equation \re{eqn_theta123}, and the $\mathrm{Sp}(n+1)$-part is given by
$\Theta=\begin{pmatrix}q^A_{~B} & t^A_{~\tilde{B}}\\-\bar{t}^{\tilde{A}}_{~B}& \bar{q}^{\tilde{A}}_{~\tilde{B}}\end{pmatrix}_{A,B=0,\ldots,n},$
where
\[ q=(q^A_{~B})_{A,B=0,\ldots,n}=\begin{pmatrix}\frac{i}{4}d^c\mathcal{K}+\frac{3}{4}(\bar{\alpha}^0-\alpha^0) & -\alpha^b\\\\\bar{\alpha}^a& \omega^a_{~b}+\frac{1}{4}(-id^c\mathcal{K}+(\bar{\alpha}^0-\alpha^0))\delta^a_{~b}\end{pmatrix}_{a,b=1,\ldots,n}\]
and
\[t=(t^A_{~\tilde{B}})_{A,B=0,\ldots,n}=\begin{pmatrix}0 & 0\\0& ie^{\mathcal{K}}\sum_{c=1}^n\tilde{h}_{abc}\,\alpha^c\end{pmatrix}_{a,b=1,\ldots,n}.\]
\label{corSpnConnQMap}
\begin{proof}
The vanishing of torsion is the following system of equations for the components of the connection one-form:
\begin{align*}
0&=d\beta^A+p^1_{~1}\wedge \beta^A-p^1_{~2}\wedge \bar{\alpha}^A+\sum_{B=0}^n(q^A_{~B}\wedge \beta^B+t^A_{~B}\wedge \alpha^B), 
\\*
0&=d\alpha^A+p^1_{~1}\wedge\alpha^A+p^1_{~2}\wedge \bar{\beta}^A+\sum_{B=0}^n(-\bar{t}^A_{~B}\wedge \beta^B+\bar{q}^A_{~B}\wedge \alpha^B) 
\end{align*}
for $A=0,\ldots,n$.
This is straightforward to solve using Proposition \ref{propositionDerivativeCoframeqMap}.
\end{proof}
\ep

\subsection{Curvature tensor for quaternionic K\"ahler manifolds obtained by the q-map}
We consider a manifold in the image of the q-map and use the notation introduced in the last section. In terms of the local frame of the type \eqref{eqQuaternionicVielbeinFrame}, the Riemann curvature tensor of a quaternionic K\"ahler manifold reads
\begin{equation*}
f^{\alpha\Gamma}(R(X,Y)f_{\beta\Delta})={\tilde{R}_H}\,\!^\alpha_{~\beta}(X,Y) \delta^\Gamma_{~\Delta}+\delta^\alpha_{~\beta}{\tilde{R}_E}\,\!^\Gamma_{~\Delta}(X,Y),
\end{equation*}
where
\begin{align*}\nonumber \tilde{R}_H&=dp+p\wedge p\\*&=\begin{pmatrix}-id\bar{\theta}_1+2i\bar{\theta}_2\wedge\bar{\theta}_3 & -(d\bar{\theta}_2+id\bar{\theta}_3)+2i\bar{\theta}_1\wedge (\bar{\theta}_2+i\bar{\theta}_3)\\(d\bar{\theta}_2-id\bar{\theta}_3)+2i\bar{\theta}_1\wedge (\bar{\theta}_2-i\bar{\theta}_3)&id\bar{\theta}_1-2i\bar{\theta}_2\wedge\bar{\theta}_3\end{pmatrix}\nonumber\\*&\stackrel{\eqref{propChQKAlekseevskyCurvature}}{=}\frac{\nu}{2}\begin{pmatrix}-i\omega_1& -\omega_2-i\omega_3\\\omega_2-i\omega_3& i\omega_1\end{pmatrix}\end{align*}
and
\be \tilde{R}_E=d\Theta+\Theta\wedge\Theta.\label{eqQKCurvE} \ee
We write the $\mathrm{Sp}(n)$-part of the curvature tensor as
\begin{equation*}
\tilde{R}_E=\begin{pmatrix}r &s\\-\bar{s} & \bar{r}
\end{pmatrix},
\end{equation*}
where $r$, $s$ are complex two-form valued $(n+1)\times (n+1)$ matrices that fulfill $r^\dagger=-r$, $s^t=s$. In terms of this splitting, equations \eqref{eqQKCOnnE} and \eqref{eqQKCurvE} read
\begin{align*}
r^A_{~B}&=dq^A_{~B}+\sum_{C=0}^n(q^A_{~C}\wedge q^C_{~B}-t^A_{~C}\wedge \bar{t}^C_{~B})
\\
s^A_{~B}&=dt^A_{~B}+\sum_{C=0}^n(q^A_{~C}\wedge t^C_{~B}+t^A_{~C}\wedge \bar{q}^C_{~B}),
\end{align*}
for $A,B=0,\ldots,n$.

Recall that in the $E$-$H$-formalism, the complexified quaternionic K\"ahler metric on the complexified tangent bundle $T^\mathbb{C}\bar{N}\cong H\otimes E$ can be written in the form $g_{FS}^\mathbb{C}=\omega_H\otimes\omega_E$, where $\omega_H$ and $\omega_E$ are non-degenerate skew-symmetric two-forms. The two-forms are represented by matrices $(\epsilon_{\alpha\beta})_{\alpha,\beta=1,2}$ and $\left(\frac{1}{2}C_{\Gamma\Delta}\right)_{\Gamma,\Delta=1,\ldots,2n+2}$, where
\begin{equation*}
	C_{\Gamma\Delta}=2\omega_E(E_\Gamma,E_\Lambda),\quad
	\epsilon_{\alpha\beta}=\omega_H(h_\alpha,h_\beta).
\end{equation*}
We have that $C_{A\tilde{B}}=-C_{\tilde{A}B}=\delta_{AB}$, $C_{AB}=C_{\tilde{A}\tilde{B}}=0$ ($A,B=0,\ldots, n$), and $\epsilon_{12}=-\epsilon_{21}=1$, $\epsilon_{11}=\epsilon_{22}=0$.


\bp
The $\mathrm{Sp}(n+1)$-part \eqref{eqQKCurvE} of the curvature tensor can be expressed as
\be {\tilde{R}_E}\,\!^{\Lambda}_{~\Xi}=\sum_{\alpha,\,\beta=1}^2\sum_{\Delta=1}^{2n+2}\frac{\nu}{4}\epsilon_{\alpha\beta}C_{\Xi\Delta}f^{\alpha\Lambda}\wedge f^{\beta\Delta}+\sum_{\alpha,\,\beta=1}^2\sum_{\Lambda',\Gamma,\Delta=1}^{2n+2}C^{\Lambda\Lambda'}\Omega_{\Lambda'\Xi\Gamma\Delta}\epsilon_{\alpha\beta}f^{\alpha\Gamma}\wedge f^{\beta\Delta},\label{equationSplittingOfCurvatureTensorInTheorem}\ee
where $C_{\Gamma\Delta}=-C^{\Gamma\Delta}$, and $\Omega_{\Lambda'\Xi\Gamma\Delta}$ are complex-valued functions on $\bar{N}$ that are symmetric in all four indices.
\begin{proof}
Recall \cite{A,S} that the curvature tensor of every $4m$-dimensional quaternionic K\"ahler manifold can be decomposed as
\begin{equation}\label{eqn_decompCT}
	R=\nu R_{\mathbb{H}P^m}+W,
\end{equation}
where
\begin{align} R_{\mathbb{H}P^m}(X,\,Y)Z=\frac{1}{4}[g(Y,\,Z)X&-g(X,\,Z)Y]-\frac{1}{2}\sum_{i=1}^3\omega_i(X,\,Y)J_i Z\nonumber\\
& +\frac{1}{4}\sum_{i=1}^3[\omega_i(Y,\,Z)J_i X-\omega_i(X,\,Z)J_i Y] \label{eq:CurvatureTensorOfHPn}
\end{align}
is the curvature tensor of the quaternionic projective space, $\nu$ is the reduced scalar curvature defined above, and $W$ is a curvature tensor of type $\mathrm{Sp}(m)$, which is related to an element in $\Omega\in\Gamma(S^4E^*)$ by the following formula:
\be W(h e,\,h' e')(h'' e'')=-\omega_H(h,\,h')\,h''\omega_E^{-1}(\Omega(e,\,e',\,e'',\cdot)), \label{eq:OmegainTheorem}\ee
$h,\,h',\,h''\in \Gamma(H),~ e,\,e',\,e''\in \Gamma(E)$. Writing the formulas \re{eq:CurvatureTensorOfHPn} and \re{eq:OmegainTheorem} in terms of our chosen frames with $m=n+1$, we obtain
\be f^{\delta\Lambda}(R^E_{\mathbb{H}P^{n+1}}(f_{\alpha\Gamma},\,f_{\beta \Delta})f_{\gamma\Xi})=-\frac{1}{4}\epsilon_{\alpha\beta}\delta^\delta_\gamma(C_{\Gamma\Xi}\delta_\Delta^{\Lambda}+C_{\Delta\Xi}\delta_\Gamma^{\Lambda}),
\label{eqCor7161}\ee
where $R^E_{\mathbb{H}P^{n+1}}$ is the $\mathrm{Sp}(n+1)$-part of $R_{\mathbb{H}P^{n+1}}$, and
\be f^{\delta\Lambda}(W(f_{\alpha\Gamma},\,f_{\beta \Delta})f_{\gamma\Xi})=-2\delta_\gamma^\delta\epsilon_{\alpha\beta}\sum_{\Lambda'=1}^{2n}\Omega_{\Gamma\Delta\Xi\Lambda'}C^{\Lambda'\Lambda}.\label{eqCor7162}\ee
Equations \re{eqCor7161} and \re{eqCor7162} now imply
\begin{equation*}
{\tilde{R}_E}\,\!^{\Lambda}_{~\Xi}(f_{\alpha\Gamma},\,f_{\beta\Delta})
=\frac{\nu}{4}\epsilon_{\alpha\beta}C_{\Xi\Delta}\delta^{\Lambda}_\Gamma-\frac{\nu}{4}\epsilon_{\beta\alpha}C_{\Xi\Gamma}\delta^{\Lambda}_\Delta-2\epsilon_{\alpha\beta}\sum_{\Lambda'=1}^{2n}\Omega_{\Gamma\Delta\Xi\Lambda'}C^{\Lambda'\Lambda}.
\end{equation*}
The above equation is equivalent to \re{equationSplittingOfCurvatureTensorInTheorem}.
\end{proof}
\ep

Using the expressions for the local Levi-Civita connection one-form given in Proposition \ref{corSpnConnQMap}, one obtains the following result (see \cite[Prop.\ 7.3.5]{D}) by inserting the $\mathrm{Sp}(n+1)$-part $\Theta$ of the Levi-Civita connection into the formula \re{eqQKCurvE} for the $\mathrm{Sp}(n+1)$-part of the curvature:
\bp\label{propWithDefOfS}
The $\mathrm{Sp}(n+1)$-part of the curvature two-form for any quaternionic K\"ahler manifold in the image of the q-map is given by $({\tilde{R}_E}\,\!^\Gamma_{~\Delta})=\begin{pmatrix}r^A_{~B} & s^A_{~\tilde{B}}\\-\bar{s}^{\tilde{A}}_{~B}& \bar{r}^{\tilde{A}}_{~\tilde{B}}\end{pmatrix}_{A,B=0,\ldots,n}$ with\footnote{All repeated lower case indices are summed over $1,\ldots,n$.}
\begin{align*}
r&=(r^A_{~B})\\
&=\begin{pmatrix}\begin{aligned}& \frac{1}{2}\big(\alpha^0\wedge \bar{\alpha}^0-\beta^0\wedge \bar{\beta}^0\\&+\sum_{C=0}^n \alpha^C\wedge \bar{\alpha}^C-\beta^C\wedge \bar{\beta}^C\big)\end{aligned} & \alpha^b\wedge \bar{\alpha}^0+\bar{\beta}^b\wedge \beta^0+ie^{\mathcal{K}}\tilde{h}_{bcd}\bar{\alpha}^c\wedge \beta^d\\&\\
\begin{aligned} \alpha^0\wedge \bar{\alpha}^a+\bar{\beta}^0\wedge \beta^a\\+ie^{\mathcal{K}}\tilde{h}_{acd}\alpha^c\wedge \bar{\beta}^d\end{aligned} & \begin{aligned} &\frac{1}{2}\delta^a_{~b}\sum_{C=0}^n(\alpha^C\wedge \bar{\alpha}^C-\beta^C\wedge\bar{\beta}^C)\\&\quad-(\beta^a\wedge\bar{\beta}^b+\bar{\alpha}^a\wedge \alpha^b)\\&\quad-e^{2\mathcal{K}}\tilde{h}_{adc}\tilde{h}_{ceb}(\alpha^d\wedge \bar{\alpha}^e+\bar{\beta}^d\wedge \beta^e)\end{aligned} \end{pmatrix}_{a,b=1,\ldots,n}
\end{align*}
and
\begin{align*}
s&=(s^A_{~\tilde{B}})\\
&=\begin{pmatrix}0 & 0\\0& ie^\mathcal{K}\tilde{h}_{abc}(\beta^0\wedge \bar{\beta}^c+\bar{\alpha}^0\wedge \alpha^c)+e^{2\mathcal{K}}\tilde{h}_{abf}\tilde{h}_{fde}\bar{\alpha}^d\wedge \beta^e-2S_{abcd}\alpha^c\wedge\bar{\beta}^d\end{pmatrix}_{a,b=1,\ldots,n},
\end{align*}
where
\begin{align}S_{abcd}&:=-\frac{1}{2}e^{2\mathcal{K}}\Big((\tilde{h}_{bcf}\tilde{h}_{fad}-4\tilde{h}_{bc}\tilde{h}_{ad})+(\tilde{h}_{acf}\tilde{h}_{fbd}-4\tilde{h}_{ac}\tilde{h}_{bd})+(\tilde{h}_{abf}\tilde{h}_{fcd}-4\tilde{h}_{ab}\tilde{h}_{cd})\nonumber\\
&\quad\quad +4\tilde{h}_a\tilde{h}_{bcd}+4\tilde{h}_b\tilde{h}_{cda}+4\tilde{h}_c\tilde{h}_{dab}+4\tilde{h}_d\tilde{h}_{abc}\Big).\label{eqS1}
\end{align}
\ep

\br{\rmfamily \normalfont
Note that the vanishing of the symmetric quartic tensor field\footnote{All repeated indices are summed over $1,\ldots,n$. Note that the symmetrization denoted by $(\ldots)$ over four indices includes a factor of $\frac{1}{4!}$.}
\begin{align}\nonumber&S_{abcd}\,\sigma^a\otimes \sigma^b\otimes \sigma^c\otimes \sigma^d\\\nonumber& =-\frac{1}{2}\frac{1}{4^3 h^2}\left(3h_{\tau(\mu\nu}\mathcal{K}^{\tau\tau'}h_{\sigma\rho)\tau'}-12h_{(\mu\nu}h_{\sigma\rho)}+16 h_{(\mu}h_{\nu\sigma\rho)}\right)dX^\mu\otimes dX^\nu \otimes dX^\sigma \otimes dX^\rho\\\nonumber
&=-\frac{1}{2}\frac{1}{4^3 h^2}\left(-12 hh_{\tau(\mu\nu}h^{\tau\tau'}h_{\sigma\rho)\tau'}-6h_{(\mu\nu}h_{\sigma\rho)}+16 h_{(\mu}h_{\nu\sigma\rho)}\right)dX^\mu\otimes dX^\nu \otimes dX^\sigma \otimes dX^\rho\\
&=: S_{\mu\nu\sigma\rho}\,dX^\mu\otimes dX^\nu \otimes dX^\sigma \otimes dX^\rho \label{eqS2}
\end{align}
on the projective special K\"ahler manifold $(\bar{M},\,g_{\bar{M}},\,J_{\bar{M}})$ is a necessary and sufficient condition for $(\bar{M},\,g_{\bar{M}})$ to be symmetric \cite{CV}.
}\er

Careful comparison of the expressions given in the above proposition with equation \eqref{equationSplittingOfCurvatureTensorInTheorem} leads to the following expression for the quartic symmetric tensor field determining the Riemann curvature tensor of a quaternionic K\"ahler manifold:
\bt\label{thQuarticTensorQMap} \!\!\!{\rm \cite[Th.\ 7.3.7]{D}}\\
For manifolds in the image of the q-map, the non-vanishing components of the quartic symmetric tensor field defined in equation \eqref{equationSplittingOfCurvatureTensorInTheorem} are given by
\[\Omega_{00\tilde{0}\tilde{0}}=\frac{1}{2},\quad\Omega_{0b\tilde{0}\tilde{d}}=\frac{1}{4}\delta_{bd},\quad\Omega_{ab\tilde{c}\tilde{d}}=\frac{1}{4}(\delta_{ac}\delta_{bd}+\delta_{ad}\delta_{bc})-\frac{1}{2}e^{2\mathcal{K}}\sum_{f=1} ^n\tilde{h}_{abf}\tilde{h}_{fcd},\]
\[\Omega_{\tilde{0}bcd}=\Omega_{0\tilde{b}\tilde{c}\tilde{d}}=-\frac{i}{2}e^{\mathcal{K}}\tilde{h}_{bcd},\quad \Omega_{abcd}=\Omega_{\tilde{a}\tilde{b}\tilde{c}\tilde{d}}=S_{abcd}\]
and symmetrization thereof, where $a,b,c,d=1,\ldots,n$.
\et

\subsection{Pointwise norm of the curvature tensor for quaternionic K\"ahler manifolds obtained by the q-map}\label{subsecNormCurv}
In this section, we give a general formula for the curvature invariant $\mathcal{S}_{W}:=\frac{1}{64}\|W\|^2\in C^\infty(\bar{N})$ for all quaternionic K\"ahler manifolds $\bar{N}=\bar{M}\times \mathbb{R}^{>0}\times \mathbb{R}^{2n+3}$ in the image of the q-map, where $W$ is the quaternionic Weyl tensor, see equation \re{eqn_decompCT}. We express $\mathcal{S}_{W}$ as a linear combination of three curvature invariants on the corresponding projective special K\"ahler manifold $\bar{M}$. Its relation to the squared pointwise norm of the Riemann curvature tensor $R$ is given by
\be \|R\|^2=80(n+1)^2+16(n+1)+64\,\mathcal{S}_{W}.\label{eqNormCurvSW}\ee
This follows from the orthogonality of the decomposition $R=\nu R_{\mathbb{H}P^{n+1}}+W$, the fact that the reduced scalar curvature is $\nu=-2$ for quaternionic Kähler manifolds obtained via the supergravity c-map, and the following formula for the squared pointwise norm of $R_{\mathbb{H}P^{n+1}}$:
\begin{equation*}
\|R_{\mathbb{H}P^{n+1}}\|^2=20n^2+44n+24=20(n+1)^2+4(n+1).
\end{equation*}
The above formula is obtained from equation \re{eq:CurvatureTensorOfHPn}.

The scalar curvature of a projective special K\"ahler manifold $\bar{M}$ in the image of the supergravity r-map is given by (see Theorem 3 in\footnote{Note that compared to \cite{CDL} we scaled the projective special K\"ahler metric $g_{\bar{M}}$ by a factor of $\frac{1}{2}$, which leads to a scaling of the scalar curvature $scal_{{\bar{M}}}$ by a factor of $2$.} \cite{CDL} in the special case $D=3$)
\begin{align}scal_{{\bar{M}}}&=-2n^2+n-2h\sum_{\alpha,\beta,\gamma=1}^n\sum_{\alpha'\!,\beta'\!,\gamma'\!=1}^nh_{\alpha\beta\gamma}h^{\alpha\alpha'}h^{\beta\beta'}h^{\gamma\gamma'}h_{\alpha'\beta'\gamma'}\nonumber\\&=-2n(n+1)+\frac{1}{32h^2}\sum_{\alpha,\beta,\gamma=1}^n\sum_{\alpha'\!,\beta'\!,\gamma'\!=1}^nh_{\alpha\beta\gamma}\mathcal{K}^{\alpha\alpha'}\mathcal{K}^{\beta\beta'}\mathcal{K}^{\gamma\gamma'}h_{\alpha'\beta'\gamma'}.\label{eqScalPVSK}\end{align}

The squared pointwise norm of the Riemann tensor of a projective special K\"ahler manifold $\bar{M}$ in the image of the r-map is
\begin{align*}
\| R_{\bar{M}}\|^2&=16 \sum_{\mu,\nu,\rho,\sigma=1}^n\sum_{\mu'\!,\nu'\!,\rho'\!,\sigma'\!=1}^nR_{\bar{\mu}\nu\sigma\bar{\rho}}\mathcal{K}^{\mu\mu'}\mathcal{K}^{\nu\nu'}\mathcal{K}^{\sigma\sigma'}\mathcal{K}^{\rho\rho'}R_{\mu'\bar{\nu}'\bar{\sigma}'\rho'},\\
&=-32\,scal_{{\bar{M}}}-32n(n+1)\\
&+\frac{1}{4^4h^4}\sum_{\mu,\nu,\sigma,\rho=1}^n\sum_{\mu'\!,\nu'\!,\sigma'\!,\rho'\!=1}^nB_{\rho\sigma\mu\nu}\mathcal{K}^{\rho\rho'}\mathcal{K}^{\sigma\sigma'}\mathcal{K}^{\mu\mu'}\mathcal{K}^{\nu\nu'}B_{\rho'\sigma'\mu'\nu'}
\end{align*}
where
\begin{equation*}
R_{\bar{\mu}\nu\sigma\bar{\rho}}=\sum_{\alpha=1}^n\mathcal{K}_{{\bar{\mu}}\alpha}R^\alpha_{~\nu\sigma\bar{\rho}}=-\mathcal{K}_{\bar{\mu}\nu}\mathcal{K}_{\sigma\bar{\rho}}-\mathcal{K}_{\bar{\mu}\sigma}\mathcal{K}_{\nu\bar{\rho}}+e^{2\mathcal{K}}\sum_{\beta,\gamma=1} ^nh_{\mu\rho\beta}\mathcal{K}^{\beta\gamma}h_{\gamma\sigma\nu}
\end{equation*}
and
\begin{equation*}
B_{\mu\nu\sigma\rho}:=\sum_{\kappa,\kappa'\!=1}^n h_{\mu\nu\kappa}\mathcal{K}^{\kappa\kappa'}h_{\kappa'\sigma\rho}.
\end{equation*}

The third real-valued function on $\bar{M}$ relevant for this discussion is
\begin{equation*}
\sum_{a,b,c,d=1}^n (S_{abcd})^2=\sum_{\mu,\nu,\sigma,\rho=1}^n\sum_{\mu'\!,\nu'\!,\sigma'\!,\rho'\!=1}^n S_{\mu\nu\sigma\rho}\mathcal{K}^{\mu\mu'}\mathcal{K}^{\nu\nu'}\mathcal{K}^{\sigma\sigma'}\mathcal{K}^{\rho\rho'}S_{\mu'\nu'\sigma'\rho'},
\end{equation*}
where the respective components are defined in equations \eqref{eqS1} and \eqref{eqS2}.

Using the quartic tensor field introduced in $\eqref{equationSplittingOfCurvatureTensorInTheorem}$, we define the following function on $\bar{N}$:
\begin{equation*}
\mathcal{S}_{W}:= \sum_{\Gamma,\Gamma'\!,\Gamma''\!,\Gamma'''\!=1}^{2n+2}\sum_{\Delta,\Delta'\!,\Delta''\!,\Delta'''\!=1}^{2n+2}\Omega_{\Gamma\Gamma'\Gamma''\Gamma'''}C^{\Gamma\Delta}C^{\Gamma'\Delta'}C^{\Gamma''\Delta''}C^{\Gamma'''\Delta'''}\Omega_{\Delta\Delta'\Delta''\Delta'''}.
\end{equation*}
Using the formulas for $\Omega$ given in Theorem \ref{thQuarticTensorQMap}, we find the following expression for $\mathcal{S}_{W}$:
\begin{align}\nonumber
\mathcal{S}_{W}&=
2\Omega_{ABCD}\Omega_{\tilde{A}\tilde{B}\tilde{C}\tilde{D}}-8\Omega_{ABC\tilde{D}}\Omega_{\tilde{A}\tilde{B}\tilde{C}D}+6\Omega_{AB\tilde{C}\tilde{D}}\Omega_{\tilde{A}\tilde{B}CD}\\
\nonumber&=2\Omega_{abcd}\Omega_{\tilde{a}\tilde{b}\tilde{c}\tilde{d}}-8\Omega_{abc\tilde{0}}\Omega_{\tilde{a}\tilde{b}\tilde{c}0}+6(\Omega_{00\tilde{0}\tilde{0}})^2+24\Omega_{0b\tilde{0}\tilde{d}}\Omega_{\tilde{0}\tilde{b}0d}+6\Omega_{ab\tilde{c}\tilde{d}}\Omega_{\tilde{a}\tilde{b}cd}\\
&\nonumber=2S_{abcd}S_{abcd}+2n(n+1)+scal_{\bar{M}}+\frac{3}{2}(n+1)+6(\frac{1}{4^3}\|R_{{\bar{M}}}\|^2+\frac{1}{4}scal_{{\bar{M}}}+\frac{n^2+n}{8})\nonumber\\&
=2S_{abcd}S_{abcd}+\frac{1}{4}(11n+6)(n+1)+\frac{3}{32}\|R_{{\bar{M}}}\|^2+\frac{5}{2}scal_{{\bar{M}}}.
\label{eqSW}\end{align}
Together with equation \eqref{eqNormCurvSW}, we obtain the following corollary:
\bc\label{Rsquarenorm}
The squared pointwise norm of the Riemann curvature tensor for any quaternionic K\"ahler manifold in the image of the q-map, defined by a cubic polynomial $h$ in $n$ variables, is
\begin{equation*}
\|R\|^2=64(n+1)(4n+3)+160\,scal_{{\bar{M}}}+6\|R_{\bar{M}}\|^2+128\sum_{a,b,c,d=1}^n(S_{abcd})^2.
\end{equation*}
\ec

\subsection{Example: A series of inhomogeneous complete quaternionic K\"ahler manifolds}
For $n\in\mathbb{N}$, we consider the following series of projective special real manifolds:
\be\mathcal{H}=\{h=1,~x>0\}\subset \mathbb{R}^{n},~h:=x\left(x^2-\sum_{i=1}^{n-1}y_i^2\right).\label{eqSeriesOfPSRMflds}\ee
Note that this corresponds to case d) of Theorem \ref{Thm2}, up to a shift in $n$. 
The coefficients $\mathcal{K}^{\bar{\mu}\nu}$ of the inverse metric of the corresponding project special K\"ahler manifold $\bar{M}$ obtained by the r-map are given by equation \re{eqn_Kinverse} in terms of the matrix
\begin{equation*}
	(h^{\mu\nu})=\frac{-1}{12x^2+4\sum\limits_{i=1}^{n-1}y_i^2}\left(\begin{matrix}
			-2x & 2y^t\\
			2y & 6x\cdot\mathbf{1}
		\end{matrix}
	\right),
\end{equation*}
where $y^t=(y_1,\ldots,y_{n-1})$.
The scalar curvature of the corresponding projective special K\"ahler manifold $\bar{M}$ in the image of the supergravity r-map can be calculated using equation \eqref{eqScalPVSK} and reads
\begin{equation*}
scal_{{\bar{M}}}=-n\cdot(2n-1)+3h\cdot \frac{n-2}{h-4x^3}+\frac{36x^3h^2}{(h-4x^3)^3}.
\end{equation*}
Furthermore, we find
\begin{align*}
\|R_{{\bar{M}}}\|^2&=\frac{16 }{\left(h-4 x^3\right)^6} \Big( h^6(n (3 n-8)+9)-4 h^5 (n (17 n-46)+57) x^3\\*&\qquad+4 h^4 (n (161 n-382)+537) x^6-64 h^3 (n (51 n-97)+99) x^9\\*&\qquad+128 h^2 (n (73 n-107)+78) x^{12}-2048 h (n (7 n-8)+3) x^{15}\\*&
\qquad+1024 n (9 n-8) x^{18}\Big)
\end{align*}
and
\begin{align*}
\sum_{a,b,c,d=1}^n (S_{abcd})^2&=\sum_{\mu,\nu,\sigma,\rho=1}^n\sum_{\mu'\!,\nu'\!,\sigma'\!,\rho'\!=1}^nS_{\mu\nu\sigma\rho}\mathcal{K}^{\mu\mu'}\mathcal{K}^{\nu\nu'}\mathcal{K}^{\sigma\sigma'}\mathcal{K}^{\rho\rho'}S_{\mu'\nu'\sigma'\rho'}\\
&=
\frac{3 x^6 }{\left(h-4 x^3\right)^6} \Big(h^4 (n (n+16)+207)-16 h^3 (n-2) (n+9) x^3\\&\qquad\qquad\qquad+96 h^2 \left(n^2+n-6\right) x^6-256 h (n-2) n x^9\\&\qquad\qquad\qquad+256 (n-2) n x^{12}\Big).
\end{align*}

\noindent
Using equation \eqref{eqSW}, the function $\mathcal{S}_{W}$ is calculated to be
\begin{align*}
\mathcal{S}_{W}
&=\frac{3}{2\left(h-4 x^3\right)^6} \Big(h^6 n (n+1)-4 h^5 (n+1) (5 n-2) x^3+8 h^4 (n (21 n+37)+112) x^6\\*&\qquad\qquad -256 h^3 (n (3 n+10)-11) x^9+256 h^2 (n (8 n+33)-20) x^{12}\\*
&\qquad\qquad -1024 h (n (3 n+11)+2) x^{15}+2048 (n+1) (n+2) x^{18}\Big)+\frac{3n}{4}(n+1).
\end{align*}
One can now check that the above function is non-constant for $n>1$. This can be seen as follows. Restricting the function to the hypersurface $(\mathbb{R}^{n}+i\mathcal{H})\times \mathbb{R}^{>0}\times \mathbb{R}^{2n+3}\subset (\mathbb{R}^{n}+iU)\times \mathbb{R}^{>0}\times \mathbb{R}^{2n+3}=\bar{N}$, we obtain a rational function of the real variable $x\geq 1$. It is now easy to check for all $n>1$ that the numerator is not proportional to the denominator.
Due to equation \eqref{eqNormCurvSW}, also the squared pointwise norm of the Riemann curvature tensor is non-constant. This shows that the quaternionic K\"ahler metrics obtained from the series of polynomials in equation \eqref{eqSeriesOfPSRMflds} are not locally-homogeneous for $n>1$. In total, we have the following:
\bt\label{nonlochomQKseries}
For $n>1$, the series of manifolds obtained from the complete projective special real manifolds in equation \eqref{eqSeriesOfPSRMflds} via the q-map consists of complete quaternionic K\"ahler manifolds that are not locally homogeneous.
\et
Note that this implies Theorem \ref{inh:Thm}, since the quaternionic K\"ahler manifolds $\bar{N}$ associated with the series \eqref{eqSeriesOfPSRMflds} admit a group of co-homogeneity one as discussed after the aforementioned theorem and in Appendix \ref{appendix_A}, see 
Example \ref{coh1Ex}.

\appendix
\section{Automorphisms of manifolds in the image of the r- and c-map}\label{appendix_A}

\subsection*{The r-map}
Let $\mathcal{H}\subset\mathbb{R}^n$ be a (connected) projective special real manifold with cubic polynomial $h$ and $\mathrm{Aut}(\mathcal{H})\subset\mathrm{GL}(n)$ its automorphism group, which consists of linear transformations that preserve the hypersurface $\mathcal{H}$. The supergravity r-map associates to $\mathcal{H}$ a projective special K\"ahler domain $\bar{M}$. This means that there exists a holomorphic function $F$ homogeneous of degree two defined on some $\mathbb{C}^*$-invariant domain $M_F\subset \mathbb{C}^{n+1}\setminus\{0\}$ such that $\bar{M}$ is the image of the Lagrangian cone
\begin{equation*}
	M=\left\{\left(z^0,\ldots,z^n,w_0,\ldots,w_n\right)^t\in M_F\times\mathbb{C}^{n+1}\subset V\ \left|\ w_I=\frac{\partial F}{\partial z^I},\ I=0,\ldots,n\right\}\right.
\end{equation*} 
under the canonical projection $V\setminus\{0\}\to P(V)$, where $V=\mathbb{C}^{2n+2}$ is endowed with the canonical symplectic structure $\sum dz^I\wedge dw_I$. As part the definition of a projective special K\"ahler domain, one does also require $\sum z^I N_{IJ} \bar{z}^J>0$ for all $z\in M_F$ and that the real symmetric matrix $\left(N_{IJ}\right):=\left(2\,\mathrm{Im}\, F_{IJ}\right)$ has signature $(1,n)$ for all $z\in M_F$. Note that such a manifold $M$ is called a \emph{conical affine special K\"ahler domain}. We define the automorphism group of $M$ as
\begin{equation*}
	\mathrm{Aut}(M):=\left.\left\{A\in\mathrm{Sp}\left(\mathbb{R}^{2n+2}\right)\subset\mathrm{Sp}\left(\mathbb{C}^{2n+2}\right)\ \right|\, AM\subset M\right\}.
\end{equation*}
The elements of $\mathrm{Aut}(M)$ preserve the affine special K\"ahler structure on $M$ induced by the embedding $M\subset V$ \cite{ACD} and, hence, also the projective special K\"ahler metric and the complex structure on $\bar{M}$. We denote by $\mathrm{Aut}(\bar{M})$ the group of holomorphic isometries of $\bar{M}$ induced by $\mathrm{Aut}(M)$. 
Recall (see Remark \ref{remHolPrepotrMap}) that for a conical affine special K\"ahler domain defined by the r-map the function $F$ takes the form $F(z^0,\ldots ,z^n)=\frac{h(z^1,\ldots ,z^n)}{z^0}$ and $M_F = \{ z_0(1,p)\mid z_0\in \bC^*, p \in \bR^n+ iU\}$, where 
$U= \bR^{>0}\cdot \mathcal{H}\subset \bR^n$ is the open cone generated by the hypersurface $\mathcal{H}$. 

Next we consider the subgroup
\begin{equation*}
	\mathrm{Aff}_{\mathcal{H}}\left(\mathbb{R}^n\right):=\left(\mathbb{R}^{>0}\times \mathrm{Aut}(\mathcal{H})\right)\ltimes\mathbb{R}^n\subset\mathrm{Aff}\left(\mathbb{R}^{n}\right)
\end{equation*}
and construct an embedding $\varphi_h:\mathrm{Aff}_{\mathcal{H}}\left(\mathbb{R}^n\right)\to\mathrm{Sp}\left(\mathbb{R}^{2n+2}\right)$ as follows. The restriction of $\varphi_h$ to the subgroup $\mathrm{Aut}(\mathcal{H})\subset\mathrm{Aff}_\mathcal{H}\left(\mathbb{R}^{n}\right)$ is defined by the canonical inclusions
\begin{equation*}
	\mathrm{Aut}(\mathcal{H})\subset\mathrm{GL}(n,\mathbb{R})\subset\mathrm{GL}(n+1,\mathbb{R})\subset \mathrm{Sp}\left(\mathbb{R}^{2n+2}\right).
\end{equation*}
Note that under these inclusions $\mathrm{GL}(n,\mathbb{R})$ acts trivially on the coordinates $z^0$ and $w_0$. When restricted to the $\mathbb{R}^{>0}$-factor, $\varphi_h$ is given by the inclusion
\begin{equation*}
	\mathbb{R}^{>0}\ni\lambda\mapsto\left(\begin{matrix}
			\lambda^{-\frac32} & 0 \\ 0 & \lambda^{-\frac12}\cdot\mathbf{1}
		\end{matrix}
	\right)\in\mathrm{GL}(n+1,\mathbb{R})\subset\mathrm{Sp}\left(\mathbb{R}^{2n+2}\right).
\end{equation*}
Finally, we define the homomorphism $\varphi_h|_{\mathbb{R}^n}:\mathbb{R}^n\to\mathrm{Sp}\left(\mathbb{R}^{2n+2}\right)$ by
\begin{equation}
	\varphi_h(v)=
		\left(\begin{array}{c|c|c|c}
				1 & 0 & 0 & 0\\ \hline
				v & \mathbf{1} & 0 & 0\\ \hline 
				-H(v,v,v) & -3H(v,v,\cdot) & 1 & -v^t\\ \hline
				3H(v,v,\cdot)^t & 6H_v & 0 & \mathbf{1}
			\end{array}
		\right),\label{eqn_varphihv}
\end{equation}
where $H\in S^3\left(\mathbb{R}^n\right)^*$ is the cubic tensor defined by $H(v,v,v)=h(v)$, $v\in\mathbb{R}^n$, and $H_v:\mathbb{R}^n\to\mathbb{R}^n$, $z\mapsto H(v,z,\cdot)^t$.

\bp
The above prescription defines an embedding
\begin{equation*}
	\varphi_h:\mathrm{Aff}_\mathcal{H}\left(\mathbb{R}^n\right)\to\mathrm{Aut}(M)\subset\mathrm{Sp}\left(\mathbb{R}^{2n+2}\right).
\end{equation*}
The induced homomorphism $\bar{\varphi}_h:\mathrm{Aff}_\mathcal{H}\left(\mathbb{R}^n\right)\to \mathrm{Aut}(\bar{M})$ is also an embedding.
\begin{proof}
It is straightforward to check that the matrix $A=\varphi_h(v)$ \re{eqn_varphihv} is symplectic, which shows that $\varphi_h$ maps into $\mathrm{Sp}\left(\mathbb{R}^{2n+2}\right)$. Similarly, one can easily verify that $\varphi_h$ is a group homomorphism. The fact that the group $\varphi_h\left(\mathrm{Aff}_\mathcal{H}\left(\mathbb{R}^n\right)\right)\subset\mathrm{Sp}\left(\mathbb{R}^{2n+2}\right)$ preserves the Lagrangian cone $M$ can be proven by checking that
\begin{equation*}
	\left.\frac{\partial F}{\partial {z}}\right|_{z'}={w}',
\end{equation*}
where $\left(\begin{matrix} z'\\w'\end{matrix}\right):=A\left(\begin{matrix} z\\w\end{matrix}\right)$.
\end{proof}
\ep
\bc Let $\mathcal{H}$ be a projective special real manifold on which $\mathrm{Aut}(\mathcal{H})$ acts with co-homogeneity $k\in \mathbb{N}_0$. 
Then the group $\bar{\varphi}_h\left(\mathrm{Aff}_\mathcal{H}\left(\mathbb{R}^n\right)\right)\subset \mathrm{Aut}(\bar{M})$ acts with co-homogeneity $k$ 
on the corresponding projective special K\"ahler domain $\bar{M}$ obtained by the r-map. 
\ec 
A similar result holds for Lie subgroups $L\subset \mathrm{Aut}(\mathcal{H})$. 
\begin{Ex} The projective special real manifolds in equation \eqref{eqSeriesOfPSRMflds} have $\mathrm{Aut}(\mathcal{H})
= \mathrm{O}(n-1)$. Thus 
$\mathrm{Aff}_\mathcal{H}\left(\mathbb{R}^n\right) \cong \left(\mathbb{R}^{>0}\times \mathrm{O}(n-1)\right)\ltimes\mathbb{R}^n$
acts with co-homogeneity one by automorphisms of the corresponding projective special K\"ahler domains $\bar{M}$ obtained by the r-map. 
\end{Ex}
\subsection*{The c-map}
Let $\bar{M}$ be a projective special K\"ahler domain of real dimension $2n$ and denote by $M\ra \bar{M}$ the corresponding 
conical affine special K\"ahler domain. The c-map associates with $\bar{M}$ a quaternionic 
K\"ahler manifold $\bar{N} = \bar{M} \times G$, where $G$ is the solvable Iwasawa subgroup 
of $\mathrm{SU}(1,n+2)$, which is of dimension $2n+4$. The quaternionic K\"ahler metric is of the form 
$g_{\bar{N}}= g_{\bar{M}} + g_G$, where $g_G$ is a family of left-invariant metrics on $G$ varying with 
$p\in \bar{M}$ \cite{CHM}. This implies the inclusion $G\subset \mathrm{Isom}(\bar{N})$. 
Moreover, the symplectic group $\mathrm{Sp}(\mathbb{R}^{2n+2})$ is a subgroup
of $\mathrm{Aut}(G)$, as can be easily seen from the structure of $G$ as solvable extension
of the $(2n+3)$-dimensional Heisenberg group. So $\mathrm{Aut}(M) \subset \mathrm{Sp}(\mathbb{R}^{2n+2})$
acts naturally on the trivial bundle $\bar{N} = \bar{M} \times G\ra \bar{M}$ mapping fibres to fibres and covering the action 
of $\mathrm{Aut}(\bar{M})$ on the base manifold. 
\bp The subgroup $\mathrm{Aut}(M)\ltimes G\subset \mathrm{Sp}(\mathbb{R}^{2n+2})\ltimes G$ acts by isometries on
$\bar{N}$. 
\begin{proof} This follows from \cite[Lemma 4]{CHM} by considering automorphisms of conical affine 
special K\"ahler domains rather than isomorphism between different 
projective special K\"ahler domains. 
\end{proof}
\ep 
\bc Let $M$ be a conical affine special K\"ahler domain and $\bar{M}$ the corresponding projective special K\"ahler domain. 
If a Lie subgroup $L\subset \mathrm{Aut}(M)$ acts with co-homogeneity $k\in \mathbb{N}_0$ on $\bar{M}$ then $L\ltimes G$ acts isometrically and with co-homogeneity $k$ 
on the corresponding quaternionic K\"ahler manifold $\bar{N}$ obtained by the c-map.
\ec 
\subsection*{The q-map}
For any quaternionic K\"ahler manifold $\bar{N}$ in the image of the q-map. We define 
 \[ \mathrm{Isom}_\mathcal{H}(\bar{N}) := \varphi_h\left(\mathrm{Aff}_\mathcal{H}\left(\mathbb{R}^n\right)\right)\ltimes G \subset \mathrm{Aut}(M) \ltimes G\subset \mathrm{Isom}(\bar{N}),\] 
 where $\mathcal{H}$ denotes the underlying projective special real manifold and $M$ the corresponding conical affine special K\"ahler domain. 
 \bc If $\mathrm{Aut}(\mathcal{H})$ acts with cohomogeneity $k\in \mathbb{N}_0$ on $\mathcal{H}$ then 
 $\mathrm{Isom}_\mathcal{H}(\bar{N})$ acts with cohomogeneity $k$ on $\bar{N}$. As a consequence, 
 $\mathrm{Isom}(\bar{N})$ has co-homogeneity $\le k$. 
 \ec 
\begin{Ex} \label{coh1Ex}Consider the quaternionic K\"ahler manifolds $\bar{N}$ associated with the projective special real manifolds in equation \eqref{eqSeriesOfPSRMflds} by the q-map. Then $\mathrm{Isom}_\mathcal{H}(\bar{N})$ acts with co-homogeneity one by isometries on $\bar{N}$. 
Note that the maximal compact subgroup of $\mathrm{Isom}_\mathcal{H}(\bar{N})$ is $\mathrm{O}(n-1)$ and that the maximal connected 
subgroup $\mathrm{Isom}_\mathcal{H}^0(\bar{N})$ has a Levi decomposition of the form 
\[ \mathrm{Isom}_\mathcal{H}^0(\bar{N}) = \mathrm{SO}(n-1) \ltimes \left( (\mathbb{R}^{>0} \ltimes \mathbb{R}^n) \ltimes G)\right), \]
where the semi-direct decomposition $(\mathbb{R}^{>0} \ltimes \mathbb{R}^n) \ltimes G$ of the radical is defined by the 
embedding $\varphi_h$. 
\end{Ex}


\begin{thebibliography}{ABCD}
\bibitem[A]{A} D.V.\ Alekseevsky, {\it Riemannian spaces with
exceptional holonomy groups}, Functional Anal.\ Appl.\ {\bf 2}
 (1968), 97--105.

\bibitem[ACM]{ACM} D.V.\ Alekseevsky, V.\ Cort\'es and T.\ Mohaupt, \emph{Conification of K\"ahler and hyper-K\"ahler manifolds}, 
Comm.\ Math.\ Phys.\ {\bf 324} (2013), no.\ 2, 637--655.

\bibitem[ACD]{ACD} D.V.\ Alekseevsky, V.\ Cort\'es and C.\ Devchand, {\it
Special complex manifolds}, J.\ Geom.\ Phys.\ {\bf 42} (2002), no.\ 1--2, 85--105.

\bibitem[ACDM]{ACDM} D.V.\ Alekseevsky, V.\ Cort\'es, M.\ Dyckmanns and T.\ Mohaupt, {\it
Quaternionic K\"ahler metrics associated with special K\"ahler manifolds}, J.\ Geom.\ Phys.\ {\bf 92} (2015), 271--287.

\bibitem[APP]{APP} S.\ Alexandrov, D.\ Persson and B.\ Pioline, \emph{Wall-crossing, Rogers dilogarithm, and the QK/HK correspondence}, JHEP1112027 (2011).

\bibitem[C]{C} V.\ Cort\'es, {\it Alekseevskian spaces}, Differential Geom.\ Appl.\ {\bf 6} (1996), no.\ 2, 129--168.

\bibitem[CDL]{CDL}
V.\ Cort\'es, M.\ Dyckmanns and D.\ Lindemann, {\it Classification of complete projective special real surfaces}, Proc.\ London Math.\ Soc.~{\bf 109}
(2014), no.\ 2, 423--445. 

\bibitem[CHM]{CHM} V.\ Cort\'es, X.\ Han and T.\ Mohaupt, {\it 
 Completeness in supergravity constructions}, Comm.\ Math.\ Phys.\ {\bf 311} 
(2012), no. 1, 191--213. 

\bibitem[CLST]{CLST}
V.\ Cort\'es, J.\ Louis, P.\ Smyth and H.\ Triendl,
{\it On certain K\"ahler quotients of quaternionic K\"ahler manifolds},
Commun.\ Math.\ Phys.\ {\bf 317} (2013), no.\ 3, 787--816.

\bibitem[CM]{CM} V.\ Cort\'es and T.\ Mohaupt, {\it Special Geometry
of Euclidean Supersymmetry III: the local r-map, instantons and
black holes}, JHEP {\bf 0907} 066 (2009).

\bibitem[CNS]{CNS}
V.\ Cort\'es, M.\ Nardmann and S.\ Suhr, {\it Completeness of hyperbolic centroaffine hypersurfaces}, Comm.\ Anal.\ Geom.\ {\bf 24} (2016), no.\ 1, 59--92.

\bibitem[CV]{CV} 
E.\ Cremmer and A.\ Van Proeyen,
{\it Classification Of K\"ahler Manifolds In N=2 Vector Multiplet Supergravity Couplings},
Class.\ Quant.\ Grav.\ {\bf 2} (1985), no.\ 4, 445--454.

\bibitem[D]{D} M.\ Dyckmanns, {\it The hyper-K\"ahler/quaternionic K\"ahler correspondence and the geometry of the c-map}, PhD thesis, University of Hamburg, 2015.

\bibitem[DS]{DS} A.\ Dancer and A.\ Swann, {\it Quaternionic K\"ahler manifolds of cohomogeneity one}, Internat.\ J.\ Math.\ {\bf 10} (1999), no.\ 5, 541--570.

\bibitem[DV]{DV} B.\ de Wit, A.\ Van Proeyen, {\it Special geometry, cubic polynomials and homogeneous quaternionic spaces}, Comm.\ Math.\ Phys.\ {\bf 149} (1992), no.\ 2, 307--333.

\bibitem[DVV]{DVV} B.\ de Wit, F.\ Vanderseypen, and A.\ Van Proeyen, {\it Symmetry structure of special geometries}, Nucl.\ Phys.\ {\bf B400} (1993) 463--524.
\bibitem[FS]{FS} S.\ Ferrara and S.\ Sabharwal, {\it 
Quaternionic manifolds for type II superstring vacua of Calabi-Yau spaces}, 
Nucl.\ Phys.\ {\bf B332} (1990), no.\ 2, 317--332.

\bibitem[F]{F} D.\ S.\ Freed, \emph{Special K\"ahler manifolds}, Comm.\ Math.\ Phys.\ {\bf 203} (1999), no.\ 1, 31--52.

\bibitem[GST]{GST} M.\ G{\"u}naydin, G.\ Sierra and P.\ K.\ Townsend, \emph{The geometry of $N=2$ Maxwell--Einstein supergravity and Jordan algebras}, Nucl. Phys. {\bf B242} (1984), 244--268.
\bibitem[Ha]{Ha} A.\ Haydys, \emph{Hyper-K\"ahler and quaternionic K\"ahler manifolds with $S^1$-symmetries}, J.\ Geom.\ Phys.\ {\bf 58} (2008), no.\ 3, 293--306.

\bibitem[Hi]{Hi} N.\ Hitchin, \emph{On the hyperk\"ahler/quaternion K\"ahler correspondence}, Commun.\ Math.\ Phys.\ {\bf 324} (2013), no.\ 1, 77--106.

\bibitem[L]{L} C.\ LeBrun, \emph{On  complete  quaternionic-K\"ahler  manifolds},  Duke  Math.\  J.\ {\bf 63} (1991), no.\ 3, 723--743.
\bibitem[MS]{MS} O.\ Mac\'ia and A.\ Swann, \emph{Twist geometry of the c-map}, Comm.\ Math.\ Phys.\ {\bf 336} (2015), 1329--1357.

\bibitem[P]{P} H.\ Pedersen, {\it Einstein metrics, spinning top motions and monopoles}, Math.\ Ann.\ {\bf 274} (1986), 35--59.

\bibitem[PV]{PV} F.\ Podest\`a, L.\ Verdiani, \textit{A note on quaternion-K\"ahler manifolds}, Internat.\ J.\ Math.\ {\bf11} (2000), no. 2, 279--283.

\bibitem[S]{S} S.\ Salamon, {\it Quaternionic K\"ahler manifolds}, Invent. Math. {\bf 67} (1982), no.\ 1, 143--171.
\end{thebibliography}
\end{document}